\documentclass[11pt]{article}
\usepackage{amsmath, upgreek}
\usepackage{amssymb}
\usepackage{amsfonts}
\usepackage{amsmath,tabu}
\usepackage{cite}
\usepackage{color}
\usepackage[final]{pdfpages}
\usepackage{xcolor} 
\usepackage{amscd}
\usepackage{xspace}
\usepackage{verbatim}
\usepackage{graphicx}
\newcommand{\mysection}[1]{
\section{#1}\setcounter{equation}{0}}
\title{\bf Local and global properties of solutions 
of an elliptic equation involving exponential and gradient reaction
}
\author{{\bf Marie-Fran\c{c}oise Bidaut-V\'eron\footnote{\noindent Laboratoire de Math\'{e}matiques et Physique Th\'{e}orique,
UMR 7013, Universit\'e de Tours, 37200 Tours, France. E-mail: veronmf@univ-tours.fr}} \\
{\bf Marta Garcia-Huidobro \footnote{\noindent
 Departamento de Matematicas, Pontifica Universidad Catolica de Chile
Casilla 307, Correo 2, Santiago de Chile. E-mail: mgarcia@mat.puc.cl}}\\
 {\bf Laurent V\'eron \footnote{\noindent
Laboratoire de Math\'{e}matiques et Physique Th\'{e}orique, UMR 7013,  Universit\'e de Tours, 37200 Tours, France. E-mail: veronl@univ-tours.fr}}\\[2mm]
}

\date{}
\begin{document}

 \maketitle


\newcommand{\txt}[1]{\;\text{ #1 }\;}
\newcommand{\tbf}{\textbf}
\newcommand{\tit}{\textit}
\newcommand{\tsc}{\textsc}
\newcommand{\trm}{\textrm}
\newcommand{\mbf}{\mathbf}
\newcommand{\mrm}{\mathrm}
\newcommand{\bsym}{\boldsymbol}
\newcommand{\scs}{\scriptstyle}
\newcommand{\sss}{\scriptscriptstyle}
\newcommand{\txts}{\textstyle}
\newcommand{\dsps}{\displaystyle}
\newcommand{\fnz}{\footnotesize}
\newcommand{\scz}{\scriptsize}
\newcommand{\be}{\begin{equation}}
\newcommand{\bel}[1]{\begin{equation}\label{#1}}
\newcommand{\ee}{\end{equation}}
\newcommand{\eqnl}[2]{\begin{equation}\label{#1}{#2}\end{equation}}
\newcommand{\barr}{\begin{eqnarray}}
\newcommand{\earr}{\end{eqnarray}}
\newcommand{\bars}{\begin{eqnarray*}}
\newcommand{\ears}{\end{eqnarray*}}
\newcommand{\nnu}{\nonumber \\}
\newtheorem{subn}{\name}
\renewcommand{\thesubn}{}
\newcommand{\bsn}[1]{\def\name{#1}\begin{subn}}
\newcommand{\esn}{\end{subn}}
\newtheorem{sub}{\name}[section]
\newcommand{\dn}[1]{\def\name{#1}}   
\newcommand{\bs}{\begin{sub}}
\newcommand{\es}{\end{sub}}
\newcommand{\bsl}[1]{\begin{sub}\label{#1}}
\newcommand{\bth}[1]{\def\name{Theorem}
\begin{sub}\label{t:#1}}
\newcommand{\blemma}[1]{\def\name{Lemma}
\begin{sub}\label{l:#1}}
\newcommand{\bcor}[1]{\def\name{Corollary}
\begin{sub}\label{c:#1}}
\newcommand{\bdef}[1]{\def\name{Definition}
\begin{sub}\label{d:#1}}
\newcommand{\bprop}[1]{\def\name{Proposition}
\begin{sub}\label{p:#1}}

\newcommand{\R}{\eqref}
\newcommand{\rth}[1]{Theorem~\ref{t:#1}}
\newcommand{\rlemma}[1]{Lemma~\ref{l:#1}}
\newcommand{\rcor}[1]{Corollary~\ref{c:#1}}
\newcommand{\rdef}[1]{Definition~\ref{d:#1}}
\newcommand{\rprop}[1]{Proposition~\ref{p:#1}}
\newcommand{\BA}{\begin{array}}
\newcommand{\EA}{\end{array}}
\newcommand{\BAN}{\renewcommand{\arraystretch}{1.2}
\setlength{\arraycolsep}{2pt}\begin{array}}
\newcommand{\BAV}[2]{\renewcommand{\arraystretch}{#1}
\setlength{\arraycolsep}{#2}\begin{array}}
\newcommand{\BSA}{\begin{subarray}}
\newcommand{\ESA}{\end{subarray}}
\newcommand{\BAL}{\begin{aligned}}
\newcommand{\EAL}{\end{aligned}}
\newcommand{\BALG}{\begin{alignat}}
\newcommand{\EALG}{\end{alignat}}
\newcommand{\BALGN}{\begin{alignat*}}
\newcommand{\EALGN}{\end{alignat*}}
\newcommand{\note}[1]{\textit{#1.}\hspace{2mm}}
\newcommand{\Proof}{\note{Proof}}
\newcommand{\qeda}{\hspace{10mm}\hfill $\square$}
\newcommand{\qed}{\\
${}$ \hfill $\square$}
\newcommand{\Remark}{\note{Remark}}
\newcommand{\modin}{$\,$\\[-4mm] \indent}
\newcommand{\forevery}{\quad \forall}
\newcommand{\set}[1]{\{#1\}}
\newcommand{\setdef}[2]{\{\,#1:\,#2\,\}}
\newcommand{\setm}[2]{\{\,#1\mid #2\,\}}
\newcommand{\mt}{\mapsto}
\newcommand{\lra}{\longrightarrow}
\newcommand{\lla}{\longleftarrow}
\newcommand{\llra}{\longleftrightarrow}
\newcommand{\Lra}{\Longrightarrow}
\newcommand{\Lla}{\Longleftarrow}
\newcommand{\Llra}{\Longleftrightarrow}
\newcommand{\warrow}{\rightharpoonup}
\newcommand{
\paran}[1]{\left (#1 \right )}
\newcommand{\sqbr}[1]{\left [#1 \right ]}
\newcommand{\curlybr}[1]{\left \{#1 \right \}}
\newcommand{\abs}[1]{\left |#1\right |}
\newcommand{\norm}[1]{\left \|#1\right \|}
\newcommand{
\paranb}[1]{\big (#1 \big )}
\newcommand{\lsqbrb}[1]{\big [#1 \big ]}
\newcommand{\lcurlybrb}[1]{\big \{#1 \big \}}
\newcommand{\absb}[1]{\big |#1\big |}
\newcommand{\normb}[1]{\big \|#1\big \|}
\newcommand{
\paranB}[1]{\Big (#1 \Big )}
\newcommand{\absB}[1]{\Big |#1\Big |}
\newcommand{\normB}[1]{\Big \|#1\Big \|}
\newcommand{\produal}[1]{\langle #1 \rangle}

\newcommand{\thkl}{\rule[-.5mm]{.3mm}{3mm}}
\newcommand{\thknorm}[1]{\thkl #1 \thkl\,}
\newcommand{\trinorm}[1]{|\!|\!| #1 |\!|\!|\,}
\newcommand{\bang}[1]{\langle #1 \rangle}
\def\angb<#1>{\langle #1 \rangle}
\newcommand{\vstrut}[1]{\rule{0mm}{#1}}
\newcommand{\rec}[1]{\frac{1}{#1}}
\newcommand{\opname}[1]{\mbox{\rm #1}\,}
\newcommand{\supp}{\opname{supp}}
\newcommand{\dist}{\opname{dist}}
\newcommand{\myfrac}[2]{{\displaystyle \frac{#1}{#2} }}
\newcommand{\myint}[2]{{\displaystyle \int_{#1}^{#2}}}
\newcommand{\mysum}[2]{{\displaystyle \sum_{#1}^{#2}}}
\newcommand {\dint}{{\displaystyle \myint\!\!\myint}}
\newcommand{\q}{\quad}
\newcommand{\qq}{\qquad}
\newcommand{\hsp}[1]{\hspace{#1mm}}
\newcommand{\vsp}[1]{\vspace{#1mm}}
\newcommand{\ity}{\infty}
\newcommand{\prt}{\partial}
\newcommand{\sms}{\setminus}
\newcommand{\ems}{\emptyset}
\newcommand{\ti}{\times}
\newcommand{\pr}{^\prime}
\newcommand{\ppr}{^{\prime\prime}}
\newcommand{\tl}{\tilde}
\newcommand{\sbs}{\subset}
\newcommand{\sbeq}{\subseteq}
\newcommand{\nind}{\noindent}
\newcommand{\ind}{\indent}
\newcommand{\ovl}{\overline}
\newcommand{\unl}{\underline}
\newcommand{\nin}{\not\in}
\newcommand{\pfrac}[2]{\genfrac{(}{)}{}{}{#1}{#2}}

\def\ga{\alpha}     \def\gb{\beta}       \def\gg{\gamma}
\def\gc{\chi}       \def\gd{\delta}      \def\ge{\epsilon}
\def\gth{\theta}                         \def\vge{\varepsilon}
\def\gf{\phi}       \def\vgf{\varphi}    \def\gh{\eta}
\def\gi{\iota}      \def\gk{\kappa}      \def\gl{\lambda}
\def\gm{\mu}        \def\gn{\nu}         \def\gp{\pi}
\def\vgp{\varpi}    \def\gr{\rho}        \def\vgr{\varrho}
\def\gs{\sigma}     \def\vgs{\varsigma}  \def\gt{\tau}
\def\gu{\upsilon}   \def\gv{\vartheta}   \def\gw{\omega}
\def\gx{\xi}        \def\gy{\psi}        \def\gz{\zeta}
\def\Gg{\Gamma}     \def\Gd{\Delta}      \def\Gf{\Phi}
\def\Gth{\Theta}
\def\Gl{\Lambda}    \def\Gs{\Sigma}      \def\Gp{\Pi}
\def\Gw{\Omega}     \def\Gx{\Xi}         \def\Gy{\Psi}

\def\CS{{\mathcal S}}   \def\CM{{\mathcal M}}   \def\CN{{\mathcal N}}
\def\CR{{\mathcal R}}   \def\CO{{\mathcal O}}   \def\CP{{\mathcal P}}
\def\CA{{\mathcal A}}   \def\CB{{\mathcal B}}   \def\CC{{\mathcal C}}
\def\CD{{\mathcal D}}   \def\CE{{\mathcal E}}   \def\CF{{\mathcal F}}
\def\CG{{\mathcal G}}   \def\CH{{\mathcal H}}   \def\CI{{\mathcal I}}
\def\CJ{{\mathcal J}}   \def\CK{{\mathcal K}}   \def\CL{{\mathcal L}}
\def\CT{{\mathcal T}}   \def\CU{{\mathcal U}}   \def\CV{{\mathcal V}}
\def\CZ{{\mathcal Z}}   \def\CX{{\mathcal X}}   \def\CY{{\mathcal Y}}
\def\CW{{\mathcal W}} \def\CQ{{\mathcal Q}}
\def\BBA {\mathbb A}   \def\BBb {\mathbb B}    \def\BBC {\mathbb C}
\def\BBD {\mathbb D}   \def\BBE {\mathbb E}    \def\BBF {\mathbb F}
\def\BBG {\mathbb G}   \def\BBH {\mathbb H}    \def\BBI {\mathbb I}
\def\BBJ {\mathbb J}   \def\BBK {\mathbb K}    \def\BBL {\mathbb L}
\def\BBM {\mathbb M}   \def\BBN {\mathbb N}    \def\BBO {\mathbb O}
\def\BBP {\mathbb P}   \def\BBR {\mathbb R}    \def\BBS {\mathbb S}
\def\BBT {\mathbb T}   \def\BBU {\mathbb U}    \def\BBV {\mathbb V}
\def\BBW {\mathbb W}   \def\BBX {\mathbb X}    \def\BBY {\mathbb Y}
\def\BBZ {\mathbb Z}

\def\GTA {\mathfrak A}   \def\GTB {\mathfrak B}    \def\GTC {\mathfrak C}
\def\GTD {\mathfrak D}   \def\GTE {\mathfrak E}    \def\GTF {\mathfrak F}
\def\GTG {\mathfrak G}   \def\GTH {\mathfrak H}    \def\GTI {\mathfrak I}
\def\GTJ {\mathfrak J}   \def\GTK {\mathfrak K}    \def\GTL {\mathfrak L}
\def\GTM {\mathfrak M}   \def\GTN {\mathfrak N}    \def\GTO {\mathfrak O}
\def\GTP {\mathfrak P}   \def\GTR {\mathfrak R}    \def\GTS {\mathfrak S}
\def\GTT {\mathfrak T}   \def\GTU {\mathfrak U}    \def\GTV {\mathfrak V}
\def\GTW {\mathfrak W}   \def\GTX {\mathfrak X}    \def\GTY {\mathfrak Y}
\def\GTZ {\mathfrak Z}   \def\GTQ {\mathfrak Q}

\font\Sym= msam10 
\def\SYM#1{\hbox{\Sym #1}}
\newcommand{\bdw}{\prt\Gw\xspace}
\maketitle\medskip

\noindent{\small {\bf Abstract} We study some local and global properties of  solutions of $-\Gd u- m\abs{\nabla u}^q-e^{u}=0$ in a punctured domain $\Gw\setminus\{0\}$, or in an exterior domain of $\BBR^N$, $N\geq 2$, where $m$ is a positive parameter  and $q>1$. We study particularly the local behaviour of solutions with an isolated singularity or the asymptotic behaviour for solutions defined in an exterior domain, and also the existence of solutions with the behaviours previously described. These behaviours change drastically according $q$ is smaller or larger than $2$. Many results are obtained by introducing various dynamical systems associated to the equation.
}\medskip

\tableofcontents
\vspace{1mm}
\hspace{.05in}
\medskip
\mysection{Introduction}
The aim of this article is to study local or global behaviour of functions which satisfy an equation of viscous {\it Chandrasekhar-Hamilton-Jacobi type} 
\bel{Na-1}
-\Gd u- m|\nabla u|^q-e^{u}=0,
\ee
either in a punctured domain of $\BBR^N$ or in an exterior domain, as well as the existence of solutions with the behaviour that we put into light. Throughout this article we assume that $m>0$, $N\geq 2$ and $q>1$. However, since our results depend strongly on the value of $q$ with respect to $\frac N{N-1}$ and $2$, in most statements we will recall the range of the parameters $q$ and $N$. Because the problem is local and invariant by translation, we assume that the equation $(\ref{Na-1})$ holds
either in $B_{r_0}\setminus\{0\}$ or in $B^c_{r_0}$. This equation exhibits a large variety of phenomena, since besides the diffusion operator,  the reaction terms are nonlinear and of a different nature. According to the range of values of the parameters and the solutions, the behaviour can be modelled by two underlying nonlinear equations:\\
the {\it viscous Hamilton-Jacobi equation} 
\bel{Na-2}
-\Gd u- m|\nabla u|^q=0,
\ee
and the {\it Emden-Chandrasekhar equation} 
\bel{Na-3}
-\Gd u-e^{u}=0.
\ee
At each occurrence the difficulty is to prove that one of the three terms in the equation is negligible compared to the other two . These two equations admit specific singular solutions which play a fundamental role 
for studying $(\ref{Na-1})$, namely for equation $(\ref{Na-3})$ solutions under the form $|x|^2e^{u(x)}=\gw(\frac{x}{|x|})$ for functions $\gw$ satisfying a nonlinear equation on $S^{N-1}$, and for equation $(\ref{Na-2})$ when $1<q<\frac{N}{N-1}$, solutions under the form $C_{N,q}|x|^{-\gb}$ with $C_{N,q}<0$ where 
\bel{beta}
\gb=\frac{2-q}{q-1}.
\ee

The following equation of a different Chandrasekhar-Hamilton-Jacobi type
\bel{Na-4}
-\Gd u+ m|\nabla u|^q-e^{u}=0
\ee
has been thoroughly studied by two of the authors in \cite{BV-V24}, \cite{BV-V25}. The main difference with $(\ref{Na-1})$ lies in the fact that in $(\ref{Na-4})$ the {\it eikonal equation } 
\bel{Na-5}
m|\nabla u|^q-e^{u}=0
\ee
had there an important role for describing the singular solutions. \smallskip

By a solution $u$ (resp. (i) supersolution, (ii) sub-solution) of $(\ref{Na-1})$ in a domain $G\subset\BBR^N$, we mean a function $u\in C^1(G)$ satisfying the following relations in the sense of distributions in $G$,
\bel{Na-1^*}
-\Gd u- m|\nabla u|^q-e^{u}=0,
\ee

\bel{Na-1^**}
\left(\!\!\!\!\!\!\!\!\!\!\!\phantom{P^{P^P}}\text {resp. }\; (i)\;-\Gd u- m|\nabla u|^q-e^{u}\geq 0\,,\,\text { resp. }\; (ii)\;- \Gd u- m|\nabla u|^q-e^{u}\leq 0.\right)\quad\quad
\ee
\smallskip

For any function expressed in spherical coordinates $u(x)=u(r,\gs)$ with $(r,\gs)\in \BBR_+\ti S^{N-1}$ we denote by $\bar u(r)$ the spherical average of $u$, and we recall that if $u$ is super-harmonic (resp. sub-harmonic) in a ball $B_R$, then $u(0)\geq \bar u(r)$ (resp. $u(0)\leq \bar u(r)$) for any $r\in (0,R)$. {\it Furthermore if $u$ is a solution of (\ref{Na-1}) (or even a supersolution) then $\bar u$ is a supersolution.} 
As a consequence of these observations, we prove

\bth{supth1} Assume $N\geq 2$ and $q>1$. \smallskip 

\nind 1- If $u$ is a supersolution of $(\ref{Na-1})$ in  $B_{r_0}\setminus\{0\}$ (resp. $B_{r_0}^c$) then $\bar u(r)$ is monotone in $B_{r_1}\setminus\{0\}$ for some $r_1<r_0$ (resp. in $B_{r_2}^c$ for some $r_2>r_0$) and it satisfies  for any $0<r'_0<r_0$ (resp. $r'_0>r_0$)
 \bel{Na-10}
r^2e^{\bar u(r)}\leq C\,\text{ in }\;B_{r'_0}\setminus\{0\}\quad\left(\text{resp }\;\text{ in }\;B^c_{r'_0}\setminus\{0\}\right),
\ee
for some $C>0$ depending on $r'_0$. \smallskip

\nind 2- If $1<q<2$, there exists no supersolution of $(\ref{Na-1})$ in an exterior domain of $\BBR^N$. If $q>2$ there exist supersolutions in any exterior domain \es

A striking difference with  $(\ref{Na-4})$ is that any supersolution of $(\ref{Na-1})$ in a punctured ball satisfies an estimate from below when $1<q<2$.
 \bth{supth1-1} Assume $N\geq 2$. If $1<q<\frac{N}{N-1}$, any supersolution $u$ of $(\ref{Na-1})$ in  $B_{r_0}\setminus\{0\}$ satisfies
 \bel{Na-10-2}
u(x)\geq -\frac 1\gb\left(\frac{1}{m(q-1)}\right)|x|^{-\gb}-\min_{|y|=\frac {r_0}2}u(y)\quad\text{if }0<|x|\leq \frac{r_0}{2}.
\ee
If $\frac{N}{N-1}\leq q<2$ any supersolution in  $B_{r_0}\setminus\{0\}$  is bounded from below in $B_{\frac{r_0}2}\setminus\{0\}$.
\es

When  supersolutions of $(\ref{Na-1})$ are bounded from below we obtain different estimates. 
 \bth{supth2} Let $q>1$ and $N\geq 2$. If $u$ is a supersolution of $(\ref{Na-1})$ in $B_{r_0}\setminus\{0\}$ locally bounded from below in $B_{r_0}$,  then $e^u$ and $|\nabla u|^q$ are locally integrable in $B_{r_0}$. There exist $c\geq 0$ and a nonnegative function $\CG\in L^1_{loc}(B_{r_0})$ such that 
\bel{Na-9x}
-\Gd u=m|\nabla u|^q+e^{u}+\CG+c_Nc\gd_0
\ee 
in $\CD'(B_{r_0})$ for some explicit constant $c_N>0$. Furthermore if $N\geq 3$ we have that  $c=0$ and $\bar u(r)=o(r^{2-N})$ as $r\to 0$, while $0\leq c<2$ and $\bar u(r)=(c+o(1))\ln\frac 1r$ if $N=2$.
 \es

Next  we set 
 \bel{kappa}
\gk=\frac{(N-1)q-N}{q-1}.
\ee
We obtain another type of estimates from below for supersolutions, and these estimates are obtained via an estimate of the gradient of $\bar u$, even if no Bernstein method can be applied as 
it will be the case with solutions. The tool is the introduction of an energy function.
\bth{supth3} Let $q>2$ and $u$ be a supersolution in $B^c_{r_0}$, then $\bar u$ satisfies
 \bel{Na-11*}
0\leq  -\bar u_r(r)\leq \left(\frac{\gk}{m}\right)^\frac{1}{q-1}r^{-\frac{1}{q-1}}\quad\text{for }r\geq r_1\,\text{ large enough},
\ee
and
 \bel{Na-11**}
\bar u(r)\geq \bar u(r_1)-\frac{q-1}{q-2}\left(\frac{\gk}{m}\right)^\frac{1}{q-1}\left(r^{\frac{q-2}{q-1}}-r_1^{\frac{q-2}{q-1}}\right)\quad\text{for all }r\geq r_1.
\ee
Finally, if $u$ is bounded from above in $B^c_{r_1}$ by some real number $\ell$, there exists $C=C(N,m,q,r_1)>0$ such that 
 \bel{Na-11***}
u(x)\geq \ell -C|x|^{\frac{q-2}{q-1}}\quad\text{for all }|x|\geq r_1.
\ee
\es

Next we give our main results concerning the solutions of $(\ref{Na-1})$. Here the Bernstein technique is an essential tool. We first prove that in the case $q>2$ {\it  isolated singularities are removable}.


\bth{supth4} Let $q>2$. If $u$ is a solution of $(\ref{Na-1})$ in $B_{r_0}\setminus\{0\}$, then there exists $C=C(N,q,m)>0$ such that 
 \bel{Na-11}
|\nabla u(x)|\leq C|x|^{-\frac 1{q-1}}\quad\text{in }B_{\frac{r_0}2}\setminus\{0\}.
\ee
As a consequence $u$ satisfies
 \bel{Na-12}
|u(x)-u(y)|\leq C'|x-y|^{-\gb}\quad\text{for all  }(x,y)\in B_{\frac{r_0}2}\setminus\{0\}.
\ee
Hence $u$ can be extended as a continuous solution of $(\ref{Na-1})$ in $B_{r_0}$.
 \es
It is important to notice that since $m>0$ and $q>2$ the gradient term has a strong regularising effect. This phenomenon points out the fact that the Hamilton-Jacobi equation is dominant and the exponential term is negligible. A similar a priori estimate of  the gradient holds for solutions in an exterior domain.

\bcor{supth4bis} Let $q>2$. If $u$ is a solution of $(\ref{Na-1})$ in $B^c_{r_0}$ there exists $C=C(N,q,m)>0$ such that 
 \bel{Na-11-bis}
|\nabla u(x)|\leq C|x|^{-\frac 1{q-1}}\quad\text{in }B^c_{2r_0},
\ee
and, with another constant $C>0$, 
 \bel{Na-12-bis}
|u(x)|\leq C|x|^{\frac {q-2}{q-1}}+\max_{|y|=r_0}|u(y)|\quad\text{in }B^c_{2r_0}.
\ee
\es
When $1<q<2$ the effect of the gradient term is not sufficient to counterbalance the exponential reaction term, thus the a priori estimate involves these two terms. This phenomenon has to be put in parallel with 
the case $m=0$ where no a priori estimate is known when $N\geq 3$.
The counterpart of \rth{supth4} when $1<q<2$ is the following,
\bth{supth5} Let $1<q<2$. If $u$ is a solution of $(\ref{Na-1})$ in $B_{\gr}(x)$, there exist $C_j=C_j(N,q,m)>0$, $j=1,2$, such that 
 \bel{Na-13}
|\nabla u(x)|\leq C_1|x|^{-\frac 1{q-1}}+C_2\max_{y\in B_{\frac{|x|}2}(x)}e^{\frac{u(y)}{2(q-1)}}\quad\text{for all }x\in B_{\frac\gr2}\setminus\{0\}.
\ee
As a consequence if $u$ is a solution in $B_{r_0}\setminus\{0\}$ and if $|x|^2e ^{u(x)}$ is locally bounded in $B_{r_0}$, the function $u$ satisfies,
 \bel{Na-14}
|u(x)|\leq \frac{C_1}{\gb}\left(|x|^{-\gb}-2^\gb r_0^{-\gb}\right)+\max_{|y|=\frac\gr 2}|u(y)|\quad\text{for all  }x\in B_{\frac{r_0}2}\setminus\{0\},
\ee
for all $\gr<r_0$.
 \es 
 
From this result and by using some very delicate techniques developed in \cite{BV-V24} we derive a precise description of the behaviour of $(\ref{Na-1})$ near an isolated singularity. For this task
we denote by $(r,\gs)\in \BBR_+\ti S^{N-1}$ the spherical coordinates in $\BBR^N$. 
\bth{supth6} Let $N\geq 3$ and $1<q<2$. If $u$ is a solution of $(\ref{Na-1})$ in $B_{r_0}\setminus\{0\}$ verifying
\bel{Na-14*}\dsps 0<\liminf_{x\to 0}|x|^2e^{u(x)}\leq \limsup_{x\to 0}|x|^2e^{u(x)}<\infty,\ee
then there exists a function $\gw\in C^2(S^{N-1})$ solution of the equation
 \bel{Na-15}
\Gd'\gw-2(N-2)+e^\gw=0\quad\text{in }S^{N-1},
\ee
where $\Gd'$ is the Laplace-Beltrami operator on $S^{N-1}$, such that 
 \bel{Na-16}
\lim_{r\to 0}\left(u(r,\gs)-2\ln\frac 1r\right)=\gw(\gs)\quad\text{uniformly in }S^{N-1}.
\ee
\es
In the assumption $(\ref{Na-14*})$ the assumption on the bound from above appears, up to now, unavoidable when $q<2$, at least when $m=0$ (see e.g. \cite {BV-V1}, \cite{Far1}). Concerning the bound from below, it has also a fundamental role since, if it is not satisfied, the behaviour of the solutions of $(\ref{Na-1})$ can be totally different.

\bth{supth6-bis} Let $N\geq 3$ and  $1<q<2$. If $u$ is a solution of $(\ref{Na-1})$ in $B_{r_0}\setminus\{0\}$ such that  $|x|^2e^{u(x)}\in L_{loc}^\infty(B_{r_0})$,
 there holds: \smallskip
 
 \nind Assertion 1- Either
  \bel{Na-16**}
\int_{B_{r_0}}|\nabla u|^qdx<\infty,
 \ee
then there exists $\gg\leq 0$ such that 
   \bel{Na-16***}
-\Gd u=m|\nabla u|^q+e^u+c_N\gg\gd_0\quad\text{in }\,\CD'(B_{r_0})
 \ee
 Furthermore 
    \bel{Na-16****}
\lim_{x\to 0}|x|^{N-2}u(x)=\gg.
 \ee
Moreover, if $\gg=0$, a property which always holds if $q\geq\frac{N}{N-1}$, then $u$ is bounded from below. Finally if $\dsps \lim_{x\to 0}|x|^2e^{u(x)}=0$ and $q\geq\frac{N}{N-1}$, the function $u$ is regular. Note thet $(\ref{Na-16**})$ holds if $q>\frac{N}{N-1}$.\smallskip

  \nind Assertion 2- Or 
  \bel{Na-17}
\int_{B_{r_0}}|\nabla u|^qdx=\infty,
 \ee
a property which can hold only if $1<q<\frac{N}{N-1}$,  then 
  \bel{Na-17*}
\liminf_{x\to 0}|x|^{N-2}u(x)=-\infty.
\ee
Moreover, if $u$ is bounded from above, the following two-side inequality is verified near $x=0$,
  \bel{Na-17**}
-\frac{1}{\gb}\left(-\frac{\gk}{ m}\right)^{\frac1{q-1}}-C_2 \leq |x|^\gb u(x)\leq -\frac{1}{\gb}\left(-\frac{\gk}{ m}\right)^{\frac1{q-1}}+C_1,
\ee
for some constants $C_j\geq 0$, j=1,2.
 \es 
 
 If $u$ is a radial solution of $(\ref{Na-1})$ in $B_{r_0}$ the assumption that $|x|^2e^{u(x)}\in L^\infty(B_{r_0})$ is unnecessary since the estimate always holds by \rth{supth1}, and then the results are more precise. The techniques we use 
 for proving the next result are based upon the representation of solutions via 3-dimensional  dynamical systems.
 
\bth{supth6-ter} Let $N\geq 3$ and $1<q<2$. If $u$ is  a radial solution of $(\ref{Na-1})$ in $B_{r_0}\setminus\{0\}$ then :\smallskip

\nind 1- Either we have
  \bel{Na-17-2}
\lim_{r\to 0}r^2e^{u(r)}=2(N-2)\quad\text{and }\;\lim_{r\to 0}r u_r(r)=-2.
\ee

\nind 2- Or $u$ is a regular solution i.e.
\bel{Na-17-3}
\lim_{r\to 0}u(r)=u_0\quad\text{for some $u_0$, and }\;\lim_{r\to 0}u_r(r)=0,
\ee
and $u_0$ can be any real number.\smallskip

\nind 3- Or $1<q<\frac N{N-1}$ and \smallskip

\nind 3-1 either 
\bel{Na-17-4}
\lim_{r\to 0}r^{\gb}u(r)=-\frac{q-1}{2-q}\left(\frac {-\gk}{m}\right)^{\frac 1{q-1}},\ee

\nind 3-2 or there exists $\gg<0$ such that 
\bel{Na-17-5}
\lim_{r\to 0}r^{N-2}u(r)=\gg.\ee
\es
 
 When $q>2$ we consider the asymptotic behaviour of solutions of $(\ref{Na-1})$ in an exterior domain under the mere assumption that $|x|^2e^{u(x)}$ is bounded.

 \bth{supth7} Let $N\geq 3$, $q>2$ and $m>0$ and  $u$ be a solution of $(\ref{Na-1})$ in $B^c_{r_0}$, then \smallskip
 
 \nind 1- If $\dsps0<\liminf_{|x|\to\infty}|x|^2e^{u(x)}\leq \limsup_{|x|\to\infty}|x|^2e^{u(x)}<\infty$, there exists a function $\gw\in C^2(S^{N-1})$ satisfying $(\ref{Na-15})$ such that 
 \bel{Na-21*} \lim_{r\to\infty}\left(u(r,\gs)+2\ln r\right)=\gw(\gs)\quad\text{uniformly in }S^{N-1}.\ee
 
 \nind 2- If $\dsps\liminf_{|x|\to\infty}|x|^2e^{u(x)}=0$, then 
 \bel{Na-22*} \lim_{r\to \infty}r^{\gb}u(r)=-\frac{q-1}{q-2}\left(\frac {\gk}{m}\right)^{\frac 1{q-1}}\quad\text{uniformly in }S^{N-1}.\ee 
\es

If $u$ is radially symmetric we obtain a similar result but with weaker assumptions.

 \bth{supth7bis} Let $N\geq 3$ and $q>2$. If $u$ is a radial solution of $(\ref{Na-1})$ in $B^c_{r_0}$, then\\
  either
 \bel{Na-21} \lim_{r\to\infty}r^2e^{u(r)}=2(N-2)\quad\text{and }\lim_{|x|\to\infty}ru_r(r)=-2,\ee 
or
 \bel{Na-22} \lim_{r\to\infty}r^{\gb}u(r)=-\frac{q-1}{q-2}\left(\frac {\gk}{m}\right)^{\frac 1{q-1}}.\ee 
\es

The existence of solutions of $(\ref{Na-1})$ with an admissible prescribed behaviour is obtained essentially in the radial cases by using construction inherited 
from the geometric theory of dynamical systems associated to hyperbolic equilibria. Concerning solutions in an exterior domain, we have the next result.

\bth{supth8} Let $N\geq 3$ and $q>2$. \smallskip
 
 \nind 1- For any $r_0>0$ there exist infinitely many radial solutions of $(\ref{Na-1})$ in $B_{r_0}^c$ satisfying  $(\ref{Na-21})$.\smallskip
 
 \nind 2- We can find $r_0>0$ and $B^*\in\BBR$ such that for any $B\leq B^*$ there exists a unique radial solution of $(\ref{Na-1})$ in $B_{r_0}^c$ satisfying $u(r_0)=B$ and 
    \bel{Na-24-bis} \lim_{r\to\infty}r^{\gb}u(r)=-\frac {q-1}{q-2}\left(\frac \gk m\right)^{\frac 1{q-1}}\quad\text{and }\,\lim_{r\to\infty}r^{\frac 1{q-1}}u_r(r)=-\left(\frac \gk m\right)^{\frac 1{q-1}}.\ee 
\es

In Appendix we give another proof of the assertion 2 with slightly different statements, clearly longer  but more general and associated to the behaviour of solutions of nonlinear dynamical systems in the neighbourhood of a saddle point equilibrium.
This proof is  based upon an extension of a construction due to Dieudonn\'e \cite{Dieu}. \\

In a similar way we prove existence of singular solutions when $1<q<2$.

 \bth{supth9} Let $N\geq 3$ and $1<q<2$. \smallskip
 
 \nind 1- There exists a unique local radial solution of $(\ref{Na-1})$ satisfying
 \bel{Na-23} \lim_{r\to 0}r^2e^{u(r)}=2(N-2)\quad\text{and }\;\lim_{|x|\to 0}ru_r(r)=-2.\ee 
 
  \nind 2- If $1<q<\frac N{N-1}$ there exists $\gr_0,\gg_0>0$ such that for $0<\gr\leq \gr_0$ and $-\gg_0<\gg<0$ there exists a radial and negative solution of $(\ref{Na-1})$ in $B_\gr\setminus\{0\}$ vanishing 
  on $\prt B_\gr$ and satisfying
   \bel{Na-24} \lim_{r\to 0}r^{N-2}u(r)=\gg\quad\text{and }\lim_{r\to 0}r^{N-1}u_r(r)=(2-N)\gg.\ee 
   
   \nind 3- If $1<q<\frac N{N-1}$, then for any $A\in\BBR$ there exists a radial solution $u_A$ of $(\ref{Na-1})$ in $B_{r_0}\setminus\{0\}$ satisfying $u_A(r_0)=A $ and 
    \bel{Na-24-bis} \lim_{r\to 0}r^{\gb}u_A(r)=-\frac{q-1}{2-q}\left(-\frac \gk m\right)^{\frac 1{q-1}}\quad\text{and }\;\lim_{r\to 0}r^{\frac 1{q-1}}u_{A\,r}(r)=\left(-\frac \gk m\right)^{\frac 1{q-1}}.\ee 

\es

When $q>2$, we also prove the existence of solution continuous but with a singularity in the gradient near $0$ as predicted in \rth{supth4}. \smallskip

\bth{supth10} Let $N\geq 2$ and $q>2$. Then for any $u_0\in\BBR$ there exists a radial solution $u$ of $(\ref{Na-21})$ defined near $0$ satisfying\smallskip
    \bel{Na-25} 
r^{\frac{2}{q-1}}(u(r)-u_0)=-\frac{q-1}{q-2}\left(\frac{\gk}{m}\right)^{\frac{1}{q-1}} +o(1)\quad\text{as }\,r\to 0.
\ee
\es

In the next section we gather the results already known in the case $q=2$ where the phenomena depends essentially on $m$.
\mysection{The case $q=2$, $N\geq 3$}
If we put $v=m^me^{mu}$ equation $(\ref{Na-1})$ with $q=2$ in the variable $u$ is transformed into a Lane-Emden equation in the positive variable $v$ with exponent $Q=\frac {1+m}{m}$,
 \bel{M1} \Gd v+v^{Q}=0,\ee
with critical exponents $Q=\frac N{N-2}$  if $m=\frac {N-2}{2}$ and $Q=\frac {N+2}{N-2}$ if $m=\frac {N-2}{4}$.
 We recall the known results obtained between 1980 and 1991 which have to be compared with our new results when $q\neq 2$.\smallskip
 
 \nind 1- If $m>\frac {N-2}2$ any positive solution $v$ of $(\ref{M1} )$ in $B_{r_0}\setminus\{0\}$ satisfies for some $\gg\geq 0$,
   \bel{M2} 
   \lim_{x\to 0}|x|^{N-2}v(x)=\gg,\quad\text{equivalently }\,  \lim_{x\to 0}\left(u(x)-\frac{N-2}{m}\ln\frac 1r\right)=k\in\BBR.
   \ee
 If $\gg=0$, then $u$ is $C^\infty$ in $B_{r_0}$ (see \cite{GiTr}).\smallskip
 
 \nind 2- If $m=\frac {N-2}2$, any positive solution $v$ of $(\ref{M1} )$ either is $C^\infty$ in $B_{r_0}$ or satisfies (see \cite{Avi})
    \bel{M3} \BA{lll}\dsps
   \lim_{x\to 0}|x|^{N-2}(-\ln |x|)^{\frac{2-N}{2}} v(x)=c(N)\quad\text{equivalently }\\
   \phantom{-----------}\dsps \lim_{x\to 0}\left(u(x)-\frac{N-2}{m}\ln\frac1{|x|}+\frac {N-2}{2m}\ln\left(-\ln |x|)\right)\right)=k.
   \EA\ee
   
   \nind 3- If $m<\frac{N-2}{2}$ there exists an explicit positive solution $v^*$ of $(\ref{M1})$ under the form  $v^*(x)=c^*(N,m)|x|^{-2m}$ and therefore an explicit solution $u^*$ 
   of $(\ref {Na-1})$ under the form $|x|^2e^{u^*}=\tilde c^*(N,m)$. \smallskip

 \nind 4- If $\frac {N-2}4<m<\frac{N-2}{2}$, any positive solution $v$ of $(\ref{M1})$ either is $C^\infty$ in $B_{r_0}$ or satisfies (see \cite{GS})
     \bel{M4} 
   \lim_{x\to 0}|x|^{2m} v(x)=c(N,m)\quad\text{equivalently }\,  \lim_{x\to 0}\left(u(x)+2\ln |x|\right)=\tilde C(N,m)(1+o(1)).
   \ee
   Furthermore there exists no global positive smooth solution in $\BBR^N$\smallskip
   
    \nind 5- If $m=\frac {N-2}4$, any positive solution $v$ of $(\ref{M1} )$ in $B_{r_0}\setminus\{0\}$ is asymptotically radial in the sense that there exists a radial solution $v^*(r)$
    of $(\ref{M1} )$  such that (see \cite{CGS})
     \bel{M5} 
v(x)=v^*(|x|) (1+o(1))\quad\text{equivalently }\,  u(x)=u^*(|x|)+o(1).
   \ee
   as $x\to 0$, where $v^*(r)=m^me^{u^*(r)}$. Furthermore the global  solutions in $\BBR^N$ are radial with respect to some point and explicit.\smallskip
   
\nind 6- If $0<m<\frac {N-2}4$, any positive solution $v$ of $(\ref{M1} )$ such that $|x|^{2m}v(x)$ (equivalently $|x|^{2}e^{u(x)}$) is bounded near $0$ satisfies 
     \bel{M6} 
   \lim_{r\to 0}r^{2m} v(r,\gs)=\gw(\gs)\,\text{ uniformly on  }\, S^{N-1}\quad\text{equivalently }\,u(x)+2\ln |x|=\ln\frac{\gw^{\frac 1m}}{m}+o(1)
   \ee
   where $\gw$ is a nonnegative solution of 
        \bel{M7} 
\Gd'\gw-2m\left(N-2-2m\right)\gw+\gw^{Q}=0\,\text{ on  }\; S^{N-1}.
   \ee
   If $\gw=0$, $u$ is $C^\infty$ in $B_{r_0}$ (see \cite{BV-V1}). \smallskip
 
 \nind It is also proved that solutions $v$ of $(\ref{M1})$ vanishing on $\prt B_{r_0}$ with all the described above behaviours do exist.\smallskip
 
 Concerning the behaviour in an exterior domain the following is known.\smallskip
 
 \nind 7- If $m\geq\frac {N-2}2$ there exists no solution (actually even supersolution) of $(\ref{Na-1})$ in an exterior domain (see \cite{Ni-Ser}).
\smallskip
 
 \nind 8- If $\frac {N-2}4<m<\frac{N-2}{2}$ any positive solution $u$ of $(\ref{Na-1})$ in $B^c_{r_0}$ satisfies $(\ref{M4})$ when $|x|\to\infty$ instead of $x\to 0$ (see \cite{GS}).
\smallskip
 
 \nind 9- If $m=\frac {N-2}4$ and $u$ is a positive solution $u$ of of $(\ref{Na-1})$  in $B^c_{r_0}$ there exists a radial solution $u^*(r)$
    of $(\ref{Na-1})$ such that (see \cite{CGS})
     \bel{M5-ext} 
u(x)=u^*(|x|) +o(1)\quad \text{as }|x|\to \infty.
   \ee
   
\nind 10- If $0<m<\frac {N-2}4$ and  $u$ is a positive solution of $(\ref{Na-1})$   in $B^c_{r_0}$ such that $|x|^{2}e^{ u(x)}$ is bounded at infinity 
there exists a positive solution $\gw$ of   $(\ref{M7})$ such that $(\ref{M6})$ holds when $|x|\to \infty$.
 
 \mysection{A priori estimates}
 \medskip

  \subsection{Local and global properties of supersolutions}
 
  \subsubsection {Proof of \rth{supth1}}
If $u$ is a supersolution then $|\overline{\nabla u}|^q\geq |\bar u_r|^q$ and $\overline{e^u}\geq e^{\bar u}$, hence $\bar u$ satisfies
               \bel{L3} 
-\bar u_{rr}-\frac{N-1}{r}\bar u_r\geq e^{\bar u}+m|\bar u_r|^q\quad \text{in }(r_0,\infty).
   \ee
Then $\bar u$ cannot have any local minimum and hence it is monotone either near $0$ if $u$ is defined in $B_{r_0}\setminus \{0\}$ or near $\infty$ if it is defined in $B_{r_0}^c$. Now by  \cite[Theorem 1-(2)]{BV-V25} there holds 
\bel{LX}e^{\bar u(r)}\leq cr^{-\min\{s,2\}},\ee 
either for $0<r\leq r'_0<r_0$ or for $r>r'_0>r_0$ according the domain where $u$ is considered and this for any $s>1$ (the sign of $m$ has no importance there). Then   in particular $r^2e^{\bar u(r)}\leq c$. Notice that if the equation is considered in 
$B_{r_0}^c$ the function $\bar u$ is decreasing and tends to $-\infty$ when $r\to\infty$. Indeed at each $r^*$ such that $\bar u_r(r^*)=0$, we have $\bar u_{rr}(r^*)<0$, and $\bar u_r(r)\to-\infty$ by $(\ref{LX})$. \\
Next we assume that $q<2$.
Following the method of \cite{BV-V25} we introduce the functions
               \bel{L4}  
x=r^2e^{\bar u}\,,\; \Phi=-r\bar u_r\,,\; t=\ln r\,,
   \ee
 and obtain the system
                  \bel{L5} \left\{\BA{lll}
x_t=x(2-\Phi)\\
\Phi_t\geq x-(N-2)\Phi+me^{(2-q)t}|\Phi|^q.
\EA\right.
   \ee
The function $x$ is bounded and $\Phi\geq 0$ satisfies
 $$\Phi_t\geq \Phi(2-N)+me^{(2-q)t}\Phi^{q}.
 $$
 If $\Phi$ is not monotone at infinity there exists a sequence $\{t_n\}$ tending to $\infty$ of local maxima $\Phi(t_n)$. From the above inequality we have $(N-2)\Phi(t_n)\geq me^{(2-q)t_n}\Phi^{q}(t_n)$, hence 
 $$\Phi^{q-1}(t_n)\leq \frac{N-2}{m}e^{(q-2)t_n}.
 $$
 This implies that $\Phi(t)\to 0$ when $t\to\infty$. Similarly, if $\Phi$ is asymptotically monotone decreasing, then $\Phi_t\leq 0$ and the previous inequality is verified for $t$ large enough and in both cases we derive that $\Phi(t)\to 0$ when $t\to\infty$. Considering the equation satisfied by $x$ we deduce that $x_t=(2+o(1))x$ and therefore $x(t)\to\infty$, which is a contradiction. Hence $\Phi$ should be asymptotically monotone increasing. Then for any $\ge>0$ we have $x-(N-2)\Phi+me^{(2-q)t}|\Phi|^q\geq (m-\ge)e^{(2-q)t}|\Phi|^q$ for $t\geq t_\ge$. Thus $\Phi$ satisfies the inequality
  $$\Phi_t\geq (m-\ge)e^{(2-q)t}\Phi^{q}\quad\text{on }\,[t_\ge,\infty).  $$
  Integrating this inequality we derive
  $$\Phi^{1-q}(t)\leq \Phi^{1-q}(t_\ge)+\frac{(m-\ge)(q-1)}{2-q}\left(e^{(2-q)t_\ge}-e^{(2-q)t}\right).
  $$
This cannot hold for any $t\geq t_\ge$ when $q<2$. Actually there would exists $T>t_\ge$ such that $\dsps\lim_{t\uparrow T}\Phi(t)=\infty$, which is again a contradiction. Hence the function $\Phi$ cannot exist for all values $t>0$.\\
Finally, when $q>2$ existence of various types of solutions of $(\ref{Na-1})$ in an exterior domain will be proved in \rth{supth8}. This ends the proof.\\
 $\phantom{-----------}$\qeda
 
 \subsubsection{Proof of \rth{supth1-1}}
 1- Let $1<q<\frac N{N-1}$. If $u$ is a supersolution of $(\ref{Na-1})$, then $U=-u$ satisfies
 \bel{LX1}
 -\Gd U+m|\nabla U|^q+e^{-U}\leq 0,
 \ee
 in the same domain. Hence the function $U_+=\max\{U,0\}$ verifies
  \bel{LX2}
 -\Gd U_++m|\nabla U_+|^q\leq 0.
 \ee
 The proof of the estimate of $U_+$, based upon Keller-Osserman method, is already known (under some variant as in \cite{Ng-Tai-Ver}) and for the sake of completeness we recall the elementary proof
of \cite[Lemma 2.3] {Ng-Tai}.  Let $0<r'_0<r_0$ and $\gm'=\max_{|y|=r'_0}U_+(y)$. Recalling that $\beta=\frac{2-q}{q-1}$, we set
\bel{L7}\psi_\ge(x)=\gl (|x|-\ge)^{-\gb}+\gm'.
\ee
If $\gl\geq \frac 1\gb\left(\frac 1{m(q-1)}\right)^{\frac 1{q-1}} $, $\psi_\ge$ is a supersolution of $(\ref{L5})$ in $B_{r'_0}\setminus B_\ge$ which dominates $U$ for $|x|=\ge$ and for $|x|=r'_0$. By the comparison principle $\psi_\ge\geq U_+$ therein and letting $\ge\to 0$ yields $(\ref{Na-10-2})$.\smallskip

\nind 2- Let  $\frac{N}{N-1}\leq q<2$, let $\phi_\ge$ be the solution  (obtained via  Schauder fixed point theorem) of 
\bel{L8}\left\{\BA{lll}
-\Gd\phi_\ge+m|\nabla\phi_\ge|^q=0\quad&\text{in }B_{r'_0}\setminus B_\ge\\
\phantom{-----,,-}\phi_\ge=\max_{|y|=\ge}U(y)&\text{in }\prt B_{\ge}\\[1mm]
\phantom{-----,,-}\phi_\ge=\gm'&\text{in }\prt B_{r'_0}.
\EA\right.\ee
By the maximum principle $u\leq \phi_\ge\leq\psi_\ge$. Letting $\ge\to 0$ the solution $\phi_\ge$ which is dominated by $\psi_0(x)=\gl |x|^{-\gb}+\gm'$ converges to a nonnegative solution $\phi$ of  
\bel{L9}\left\{\BA{lll}
-\Gd\phi+m|\nabla\phi|^q=0\quad&\text{in }B_{r'_0}\setminus \{0\}\\
\phantom{------,,}\phi=\gm'&\text{in }\prt B_{r'_0}.
\EA\right.\ee
By \cite[Theorem A2]{Ng-Tai-Ver} the function $\phi$ is constant. This implies the claim.\qeda\medskip\\

\subsubsection{ Proof of \rth{supth2}} 
If $u$ is bounded from below, then $u_\gm= u+\gm$ is positive for some $\gm$. Since 
$$-\Gd u_\gm\geq m|\nabla u_\gm|^q+e^{-\gm}e^{u_\gm},
$$
it follows from Brezis-Lions lemma \cite{Br-Li} that $|\nabla u_\gm|^q$ and $e^{-\gm}e^{u_\gm}$ are locally integrable in $B_{r_0}$ and there exist $c\geq 0$ and a nonnegative function 
$\CG\in L^1_{loc}(B_{r_0})$ such that 
$$-\Gd u_\gm= m|\nabla u_\gm|^q+e^{-\gm}e^{u_\gm}+\CG+c_Nc\gd_0\quad\text{in }\; \CD'(B_{{r_0}}).
$$
When $N\geq 3$ we have that $\bar u_\gm(r)= cr^{2-N}(1+o(1))$ . The integrability of $e^{u_\gm}$ implies $c=0$ and $\bar u(r)=o(r^{2-N})$ when $r\to 0$. When $N=2$ we have that 
$\bar u_\gm(r)= (c+o(1))\ln \frac 1r$ and by comparison with $E(x)=c\ln \frac {r_0}{|x|}$ which is the solution of 
$$-\Gd E=c\gd_0\quad\text{in }\; \CD'(B_{{r_0}})
$$
vanishing on $\prt B_{r_0}$ we have that $u_\gm\geq E$. The integrability of $e^{u_\gm}$ implies $c<2$. \qeda\medskip

By introducing some energy functions we obtain an estimate of the radial gradient of a supersolution when $q\geq \frac N{N-1}$.

\bprop{corth2} Let $N\geq 2$, $q\geq\frac N{N-1}$ and $u$ be a supersolution of $(\ref{Na-1})$ in $B_{r_0}\setminus\{0\}$. Then $\bar u$ is locally bounded from below 
 in $B_{r_0}\setminus\{0\}$ and $\bar u_r$ satisfies
\bel{Na-9-1}0\leq-\bar u_r(r)\leq \left(\frac\gk m\right)^{\frac 1{q-1}}r^{-\frac{1}{q-1}}(1+o(1)))\quad\text{as }r\to 0,
\ee
if $q>\frac N{N-1}$, and 
\bel{Na-9-2}
0\leq -\bar u_r(r)\leq \left(\frac{N-1}{m}\right)^{N-1}r^{1-N}\left(-\ln r\right)^{1-N}(1+o(1)),
\ee
if $q=\frac N{N-1}$. As a consequence $\bar u(r)$  admits a limit when $r\to 0$ if $q>2$.
\es
\Proof Since
$$-\left(r^{N-1}\bar u_r\right)_r\geq mr^{-(N-1)(q-1)}|r^{N-1}\bar u_r|^q,
$$
we set $w=r^{N-1}\bar u_r$, and obtain if $q\neq \frac N{N-1}$
\bel{L1-1}
\frac{d}{dr} G(r)\geq 0\,\text{ with }\;G(r):=|w|^{-q}w+\frac{m}{\gk}r^{N-q(N-1)}.
\ee
Assume first $q> \frac N{N-1}$, then  
\bel{L1-2}
|w|^{-q}w(r)+\frac{m}{\gk}r^{N-q(N-1)}\leq G(r_0)\quad\text{for }0<r\leq r_0.
\ee
Hence
$$0\leq -w(r)\leq\left(\frac\gk m\right)^{\frac 1{q-1}}r^{\frac{q(N-1)-N}{q-1}}(1+o(1))),
$$
which implies 
$$0\leq-\bar u_r(r)\leq \left(\frac\gk m\right)^{\frac 1{q-1}}r^{-\frac{1}{q-1}}(1+o(1)))
$$
for $0<r\leq r_1$. Integrating this inequality on $(r,s)$ for $0<r<s\leq r_1$, we derive
\bel{L1-3}
0\leq \bar u(r)-\bar u(s)\leq C\left(s^{\frac{q-2}{q-1}}-r^{\frac{q-2}{q-1}}\right)\quad\text{for }0<r\leq s,
\ee
provided $q\neq 2$. Therefore, if $q>2$, $\bar u(r)$ admits a limit when $r\to 0$.\\
If $q=\frac N{N-1}$, $(\ref{L1-1})$ is replaced by
\bel{L1-5}
\frac{d}{dr} G^*(r)\geq 0\,\text{ with }\;G^*(r):=|w|^{-\frac N{N-1}}w+\frac m{N-1}\ln\frac 1r.
\ee
Hence $w(r)<0$ for $0<r\leq r_1\leq\min \{1,r_0\}$. This implies $(\ref{Na-9-2})$
and leads us to
\bel{L1-7}
\bar u(r_1)\leq \bar u(r)\leq \left(\frac{N-1}{m}\right)^{N-1}\frac{r^{2-N}}{N-2}\left(-\ln r \right)^{1-N}(1+o(1)).
\ee
This ends the proof.
\qeda

\subsubsection { Proof of \rth{supth3}}
Let $u$ be a supersolution in $B_{r_0}^c$. By $(\ref{Na-10})$ there holds
$$\bar u(r)\leq -2\ln r+\ln C\quad\text{for }r>2r_0.
$$
 Then $\bar u(r)\to-\infty$ and $\bar U(r)=-\bar u(r)\to\infty$ when $r\to\infty.$  Hence $\bar U_r\geq 0$ for $r\geq r_1<r_0$. We set $W=r^{N-1}\bar U_r$, then 
$-W_r+mr^{(N-1)(1-q)}W^q\leq 0$ which implies that the function
\bel{L15}r\mapsto H(r)=W^{1-q}-\frac{m}{\gk}r^{N-(N-1)q}
\ee
is decreasing on $(r_0,\infty)$. Since $q>2$, we have that $\gk>0$ and
$$\liminf_{r\to\infty}H(r)=\lim_{r\to\infty}H(r)=\lim_{r\to\infty}W^{1-q}(r)=\gl\geq 0.$$
If $\gl>0$ this implies that $\bar U_r$ is integrable on $(r_0,\infty)$, hence there exists a finite limit to $\bar U(r)$ when $r\to\infty$, contradiction. Hence $H(r)\to 0$ when $r\to\infty$ and thus
it is positive. Therefore there holds
$$W^{1-q}(r)\geq \frac{m}{\gk}r^{N-(N-1)q}
$$
for $r\geq r_0$, equivalently
$$-\bar u_r(r)=\bar U_r(r)\leq \left(\frac\gk m\right)^{\frac 1{q-1}}r^{-\frac{1}{q-1}}.
$$
By integration we obtain $(\ref{Na-11**})$.\\
Finally suppose that $u$ is bounded from above in $B_{r_1}^c$ by some constant $\ell$ and set $U_\ell=\ell -u$. This function is nonnegative and subharmonic thus we can control its pointwise value from above by its spherical average. From \cite[Lemma 2.1]{BV-Gr}, and since $U_\ell$ is a positive subharmonic function we can control from above the pointwise value of $U_\ell$ by its spherical average $\bar U_\ell$ in the sense that 
\bel{L17}\bar U_\ell(x)\leq C_\ge\max \{\bar U_\ell((1+\ge)|x|),\bar U_\ell((1-\ge)|x|)\}
\ee
for $0<\ge<\ge_0$ and $C_\ge>0$. This implies 
\bel{L18}\BA{lll}
 u(x)\geq \ell-C|x|^{\frac{q-2}{q-1}}\quad\text{for all }|x|\geq r_1.
\EA\ee
\qeda

 \subsection{Pointwise estimates of solutions}
 
\subsubsection{ Proof of \rth{supth4}, \rcor{supth4bis} and \rth{supth5}}
 {\it Proof of \rth{supth4}}. The proof is based upon Bernstein technique. We set $z=|\nabla u|^2$ and we have classically
 $$-\frac 12\Gd z+\frac{(\Gd u)^2}{N}+\langle\nabla (\Gd u),\nabla u\rangle \leq 0.
 $$
 Using $(\ref{Na-1})$ we obtain
\bel{P1}\BA{lll}\dsps
 -\frac 12\Gd z+\frac{(m z^{\frac q2}+e^u)^2}{N}\leq e^uz+\frac{mq}{2}z^{\frac q2-1}\langle\nabla z,\nabla u\rangle\\[0mm]
 \phantom{\dsps
 -\frac 12\Gd z+\frac{(m z^{\frac q2}+e^u)^2}{N}}\dsps
 \leq e^uz+\frac{m^2}{2N}z^q+\frac{Nq^2}{2}\frac{|\nabla z|^2}{z},
\EA \ee
since $z^{\frac {q-1}2}|\nabla z|=\frac{|\nabla z|}{\sqrt z}z^{\frac{q}2}\leq \frac{m}{Nq} z^q+\frac{Nq}{m}\frac{|\nabla z|^2}{z}$. Next
$$\frac{2mz^{\frac q2}e^u}{N}+\frac{e^{2u}}{N}-e^uz\geq \frac{e^{2u}}{N}-\frac{q-2}{q}\left(\frac{N}{mq}\right)^{\frac{2}{q-2}}e^u\geq-\frac{N(q-2)^2}{4q^2}\left(\frac{N}{mq}\right)^{\frac{4}{q-2}} :=-C_{N,q,m}.
$$
Hence
\bel{P2} -\frac 12\Gd z+\frac{m^2}{2N}z^q\leq C_{N,q,m}+\frac{Nq^2}{2}\frac{|\nabla z|^2}{z}.
\ee
Using \cite [Lemma 3.1]{BV-V2} (see also \cite[Lemma 2.2, Theorem A]{BVGHV}) we obtain that 
\bel{P3}
z(x)\leq c_1|x|^{-\frac 2{q-1}}+c_2\quad\text{for }0<|x|\leq\frac {r_0}2,
\ee
for some $c_1,c_2$ depending on $N,q,m$. This implies $(\ref{Na-11})$ and inequality $(\ref{Na-12})$ follows by integration.\qeda
\medskip

\nind{\Remark} Estimate $(\ref{Na-11})$ has the remarkable feature that the a priori estimate holds for any solution, even in presence of the exponential forcing term. A property which cannot hold if $m=0$ (see \cite{BV-V1}). This means that the gradient reaction term has a regularising effect on the exponential reaction. Moreover, this phenomenon does not depend on the sign of $m$ since it was already observed in \cite{BV-V24} for equation $(\ref{Na-4})$
with $m>0$ and $q>2$.\medskip

\nind{\it Proof of \rcor{supth4bis}} Let $x_0\in B^c_{2r_0}$, then $B_{|x_0|-r_0}(x_0)\subset B^c_{r_0}$. We apply estimate $(\ref{Na-11})$ in $B_{|x_0|-r_0}(x_0)$ and we have for any 
$x\in B_{\frac {|x_0|-r_0}2}(x_0)$,
$$|\nabla u(x)|\leq C|x-x_0|^{-\frac1{q-1}},
$$
where $C=C(N,m,q)$. If we take in particular $|x_0-x|=|x|-r_0=\frac {|x_0|-r_0}{2}$, we obtain.
\bel{P2*}
|\nabla u(x)|\leq C(|x]-r_0)^{-\frac1{q-1}}\leq 3^{\frac 1{q-1}}C|x|^{-\frac1{q-1}},
\ee
and this holds for any $x\in B_{\frac{3r_0}2}^c$. Estimate $(\ref{Na-12})$ follows by integration.\qeda

\medskip

\nind{\it Proof of \rth{supth5}}. The proof is a variant of the one of \rth{supth4}. Inequality $(\ref{P1})$ is still valid, but the term $e^uz$ on the right-hand side cannot be absorbed and $(\ref{P1})$ is replaced by 
\bel{P4}\BA{lll}\dsps -\frac 12\Gd z+\frac{m^2}{4N}z^q\leq \frac{q-1}{q}\left(\frac{4N}{qm^2}\right)^{\frac1{q-1}}e^{\frac{qu}{q-1}}+\frac{Nq^2}{2}\frac{|\nabla z|^2}{z}\\[4mm]
\phantom{-\frac 12\Gd z+\frac{m^2}{4N}z^q}\dsps
\leq \frac{q-1}{q}\left(\frac{4N}{qm^2}\right)^{\frac1{q-1}}\!\!\!\!\!\!\!\max_{|x-y|\leq \frac {|x|}2}e^{\frac{qu(y)}{q-1}}+\frac{Nq^2}{2}\frac{|\nabla z|^2}{z}
\EA\ee
Using again \cite [Lemma 3.1]{BV-V2} we obtain $(\ref{Na-13})$. If we assume that $|x|^2e^{u(x)}$ remains bounded near $0$, $(\ref{Na-14})$ follows by integration.\qeda
\mysection{Local behaviour of solutions}

The proofs of some of the convergence theorems are adaptations to equation $(\ref{Na-1})$ of results concerning singular solutions of $(\ref{Na-4})$. We recall some intermediate results the proof of them is given in \cite{BV-V24}.  We denote by $(r,\gs)\in \BBR_+\ti S^{N-1}$ the spherical coordinates in $\BBR^N$. If $u$ is a solution of $(\ref {Na-1})$  in $B_{r_0}\setminus\{0\}$ it satisfies
\bel{Q1}u_{rr}+\frac {N-1}{r}u_r+\frac1{r^2}\Gd'u+(e^{u}+m\left(u_r^2+r^{-2}|\nabla 'u|^2\right)^{\frac q2}=0\quad\text{in }(0,r_0)\ti S^{N-1}.
\ee
We set
\bel{Q2}r=e^{-t}\Longleftrightarrow t=\ln\frac 1{r}\quad\text{and }\,v(t,\gs)=u(r,\gs)-2\ln\frac 1{r}\Longrightarrow v(t,\gs)=u(r,\gs)-2t,\ee
and there holds in $(t_0,\infty)\ti S^{N-1}$ (with $t_0=-\ln r_0$)
\bel{Q3}v_{tt}-(N-2)v_t+\Gd'v-2(N-2)+e^{v}+me^{(q-2)t}\left(\left(v_t+2\right)^2+|\nabla 'v|^2\right)^{\frac q2}=0.
\ee
If $u$ is a solution of $(\ref{Na-1})$ in $B_{r_0}\setminus\{0\}$ such that $|x|^2e^u\leq\tilde K_1$, then the function $v$ satisfies 
\bel{Q4} v(t,\gs)\leq  K_1.
\ee

The next regularity result is proved in \cite[Lemma 3.2]{BV-V24} is valid for any real number $m$.
\blemma{Reg} Assume $1<q<2$, $m$ is any real number and $u$ is a solution $(\ref{Na-1})$ in $B_{r_0}\setminus\{0\}$ such that the function $v$ defined by $(\ref{Q2})$ satisfies 
\bel{Q5} 
-\eta(t)\leq v(t,\gs)\leq K_1\quad\text{for }t\geq t_0,
\ee
for some positive function $\eta\in C^2([t_0,\infty))$ verifying 
\bel{Q6} 
\frac{1}{\eta(t)}+\left|\frac{\eta_{t}(t)}{\eta(t)}\right|+\left|\frac{\eta_{tt}(t)}{\eta(t)}\right|\leq C\quad\text{and }\;\eta(t)\leq Ce^{\gb t}\quad\text{for }t\geq t_0,
\ee
then 
\bel{Q7} 
\left|v_t(t,\gs)\right|+\left|\nabla 'v(t,\gs)\right|\leq C\eta(t)\quad\text{for }(t,\gs)\in [t_1,\infty)\ti S^{N-1}.
\ee
\es
The next isotropy estimate proved in \cite[Proposition 3.3]{BV-V24} gives an estimate on $u(r,.)-\bar u(r)$ where $\bar u(r)$ is the spherical average of $u(r,.)$. Its proof combines 
the integral representation of $u-\bar u$ via semigroup associated to fractional operators, regularising effect from $L^2(S^{N-1})$ into  $L^\infty(S^{N-1})$ and some Fourier analysis.

\blemma{isotropy} Let $N\geq 3$, $1<q<2$ and $u$ be a solution $(\ref{Na-1})$ in $B_{r_0}\setminus\{0\}$ such that the function $v$ defined by $(\ref{Q2})$ satisfies 
\bel{Q8} 
-Ce^{\gn t}\leq v(t,\gs)\leq K_1\quad\text{for }t\geq t_0,
\ee
for some $0\leq \gn\leq \gb$ satisfying also $\gn<\min\left\{N-1,\frac{N+1}{q}-1\right\}$, then 
\bel{Q9} \dsps
|u(x)-u(y)|\leq K_3r^{2-q-q\gn}+K_4\quad\text{for all }0<|x|=|y|=r\leq\frac {r_0}2, 
\ee
for some $K_3,K_4>0$.
\es
\medskip

\subsection{ Characterisation of isolated singularities}

 We assume that $1<q<2$ and $u$ is a solution of $(\ref{Na-1})$ such that $|x|^2e^{u(x)}\in L^\infty(B_{r_0})$.\smallskip

\subsubsection{ Proof of \rth{supth6}}
If we assume that $\dsps 0<\liminf_{x\to 0}|x|^2e^{u(x)}\leq \limsup_{x\to 0}|x|^2e^{u(x)}<\infty$, the function $u(x)-2\ln\frac 1{|x|}$ is bounded in $B_{\frac{r_0}{2}}$. Hence the function $v$ defined in 
$(\ref{Q2})$ is bounded. Therefore the functions $v_t,\nabla 'v, v_{tt}, \nabla 'v_t$ and $\nabla'^2v$ (here $\nabla'^2v$ is the spherical Hessian on $S^{N-1}$ of the function $v(t,.)$) are uniformly bounded in $[t_1,\infty)\ti S^{N-1}$ by standard regularity results. Furthermore they are uniformly continuous there. 
The following energy identity holds:
\bel{Q10} \BA{lll}\dsps
\frac{d}{dt}\int_{S^{N-1}}\left(\frac{1}{2}(v_t^2-|\nabla 'v|^2)-2(N-2)+e^v\right) dS\\[4mm]
\phantom{----}\dsps=(N-2)\int_{S^{N-1}}v_t^2dS-me^{(q-2)t}\int_{S^{N-1}}((v_t+2)^2+|\nabla 'v|^2)^{\frac q2}v_tdS.
\EA\ee
The bounds on $v$ and its derivatives imply that $v_t\in L^2((t_1,\infty)\ti S^{N-1})$. By differentiating the equation satisfied by $v$ and using the above result it follows that 
$v_{tt}\in L^2((t_1,\infty)\ti S^{N-1})$, which combined with the uniform continuity of $v_t$ and $v_{tt}$ implies that 
\bel{Q11} \dsps
\lim_{t\to\infty}(\norm {v_t}_{L^2(S^{N-1})}+\norm {v_{tt}}_{L^2(S^{N-1})})=0.
\ee
The trajectory of $v$ in $C^2(S^{N-1})$ is defined by $\dsps\CT[v]:=\cup_{t>0}\{v(t,.)\}$. Its limit set expressed by
$$\dsps \Gw[v]=\bigcap_{t>0}cl_{C^2}\left(\bigcup_{\gt\geq t}\{v(t,.)\}\right),
$$
is a non-empty compact and connected subset of the set of solutions $\gw$ of $(\ref{Na-15})$. Using as in \cite{BV-V24}
the Huang-Takac extension of Simon's result \cite{HuTak} of asymptotics of analytic functionals  we obtain that $\Gw[v]=\{\gw\}$ for some $\gw$ satisfying $(\ref{Na-15})$, which implies the result.\qeda.

\subsubsection{ Proof of \rth{supth6-bis}} 
We first prove the following variant of \cite[Proposition 2.7]{BV-V24}.
\blemma {integ} Let $N\geq 3$ and $1<q<2$. Let $u\in C^2(B_{r_0}\setminus \{0\})$ be a solution of $(\ref{Na-1})$ such that $|x|^2e^u\in L^\infty(B_{r_0})$. If  $|\nabla u|^q\in L^1(B_{r_0})$, then $e^u$ is also integrable in $B_{r_0}$ and there exists 
$\gg\leq 0$ such that 
\bel{Q12}
-\Gd u= m|\nabla u|^q+e^u+c_N\gg\gd_0\quad\text{in }\CD'(B_{r_0}).
\ee
Furthermore,  there holds
\bel{Q12*}\BA{lll}\dsps
\min_{\prt B_{r_0}}u+\gg (|x|^{2-N}-r_0^{2-N}) \leq  u(x)\\[2mm]\phantom{-----}
\dsps\leq\min\left\{ 2\ln\frac1{ |x|}+K_1,\gg (|x|^{2-N}-r_0^{2-N})+\BBG_{B_{r_0}}[e^u+m|\nabla u|^q](x)\right\}
\EA\ee
in $B_{r_0}$, where $\BBG_{B_{r_0}}[ \; \!.  \;\!]$ is the Green operator in $B_{B_{r_0}}$.  If $q>\frac N{N-1}$ we have always $|\nabla u|^q\in L^1(B_{r_0})$.
\es
\Proof By assumption $u(x)\leq 2\ln\frac1{ |x|}+K_1$ for some $K_1>0$. We set  
$$V(x)=-u(x)+\ln\frac 1{|x|^2}+K_1,
$$
then $V\geq 0$
\bel{S7}
-\Gd V-\frac{2(N-2)}{|x|^2}+\frac{1}{|x|^2}e^{K_1-V}+m\left|\nabla u\right|^q=0.
\ee
Since $N\geq 3$ and $V\geq 0$, both $\frac{2(N-2)}{|x|^2}$ and $\frac{1}{|x|^2}e^{K_1-V}$ are integrable. By the Brezis-Lions Lemma since $\left|\nabla u\right|^q$ and the terms involving $|x|^{-2}$ are integrable in $B_1$, there exists $\tilde \gg\geq 0$ such that 
\bel{S8}
-\Gd V-\frac{2(N-2)}{|x|^2}+\frac{1}{|x|^2}e^{K_1-V}=m\left|\nabla u\right|^q+c_N\tilde \gg\gd_0\quad\text{in }\CD'(B_1).
\ee
Replacing $V$ by $-u(x)+\ln\frac 1{|x|^2}+K_1$ it follows that  $(\ref{Q12})$ holds with $\gg=-\tilde\gg$.  The proof of $(\ref{Q12*})$ is similar as the one of \cite[Proposition 2.7]{BV-V24}. If
$\gg<0$ we have $\bar u(|x|)\sim \gg |x|^{2-N}$.
Finally, if $q>\frac N{N-1}$, it follows from \rth{supth5} that $|\nabla u|^q\in L^1(B_{r_0})$ .
\qeda\medskip

\nind{\it Proof of \rth{supth6-bis}: Assertion 1}. Let $(\ref{Na-16**})$ be verified. Then $e^u+m|\nabla u|^q$ is integrable and positive, and consequently $(\ref{Q12})$ and $(\ref{Q12*})$ hold for some 
$\gg\leq 0$. By averaging this last inequality, be obtain
\bel{Q12**}\BA{lll}\dsps
\min_{\prt B_{r_0}}u+\gg (|x|^{2-N}-r_0^{2-N}) \leq  \bar u(|x|)\\[2mm]\phantom{-----}
\dsps\leq\min\left\{ 2\ln\frac1{ |x|}+K_1,\gg (|x|^{2-N}-r_0^{2-N})+\overline {\BBG_{B_{r_0}}[e^u+m|\nabla u|^q](|x|)}\right\}.
\EA\ee
Since it is easy to verify that  
$$\lim_{r\to 0}r^{N-2}\overline{G_{B_{r_0}}[e^u+m|\nabla u|^q]}(r)=0,
$$
we derive that 
\bel{Q12***}\bar u(|x|)=\gg|x|^{2-N}+o(|x|^{2-N}) \quad\text{as }\;x\to 0.
\ee
Next we set $U=-u$ and get (with $\tilde\gg=-\gg$ as in \rlemma {integ})
$$-\Gd U+m|\nabla U|^q+e^{-U}=c_N\tilde\gg\gd_0\quad\text{in }\CD'(B_{r_0}).
$$
Then $U_+$ satisfies
$$-\Gd U_++m|\nabla U_+|^q\leq 0\quad\text{in }B_{r_0}\setminus\{0\}.
$$
If $q\geq \frac{N}{N-1}$ it follows, as in the proof of \rth{supth1-1}-Step 2, that $U_+$ remains bounded near $0$ and the same holds with $\bar U_+$. {\it This implies that 
$\gg=0$ and $u$ is bounded from below}. Notice that estimate holds whenever $\gg=0$.\\
If $1<q<\frac{N}{N-1}$ and $\gg<0$. We can apply the isotropy estimate $(\ref{Q9})$ with $\gn=N-2$ since $N-2<\min\{\gb,N-1,\frac{N+1}{q}-1\}$ and obtain that 
$$|x|^{N-2}|u(x)-\bar u(|x|)|\leq K_3|x|^{N-q(N-1)}+K_4|x|^{N-2}|.
$$
Combining this inequality with $(\ref{Q12***})$ implies $(\ref{Na-16***})$. We also deduce by standard regularity estimates that 
\bel{S13}
\lim_{x\to 0}|x|^{N-1} |\nabla u(x)|=\gg(N-1).
\ee 
\smallskip

\nind{\it We claim that if $|x|^2e^{u(x)}\to 0$ when $x\to 0$, then $u$ is regular}. If $K\leq u(x)\leq -2\ln |x|+K_1$ it follows by \rlemma {Reg} that 
\bel{S14}|\nabla u(x)| \leq C|x|^{-1}\ln \frac 1{|x|} \quad\text{for }0<|x|\leq \frac 12\min\{1,r_0\}=\gr_0.
\ee
Since $1<q<2$, we have that $|\nabla u|^q\in L^p(B_{r_0})$ for some $p$ such that $\frac N2<p<\frac Nq$. Thus $\BBG_{B_{r_0}}[|\nabla u|^q]\in L^\infty(B_{r_0})$. Averaging the equation  we obtain that
 $$-\left(r^{N-1}\bar u_r(r)\right)_r=k(r)r^{N-1}e^{\bar u(r)}+mr^{N-1}\overline{|\nabla u|\!\!\!\phantom{I}^q}(r),
 $$
 where the function $k$ satisfies $e^{-K'_5}\leq k(r)\leq e^{K'_5}$. When $r\to 0$, we use the fact that  $r^{N-1}\bar u_r(r)\to 0$ to derive from the previous identity that  there holds
 $$-\bar u_r(r)=r^{1-N}\int_0^rk(r)s^{N-1}e^{\bar u(s)}ds +m r^{1-N}\int_0^rs^{N-1}\overline{|\nabla u|\!\!\!\phantom{I}^q}(s)ds.
 $$
 For any $\ge>0$ there exists $r_\ge\in (0,\gr_0]$ such that for any $r\in (0,r_\ge]$ we have that $s^2e^{\bar u(s)}\leq \ge$, thus
 $$r^{1-N}\int_0^rk(r)s^{N-1}e^{\bar u(s)}ds =r^{1-N}\int_0^rk(r)s^{N-3}s^2e^{\bar u(s)}ds\leq \frac{e^{K'_5}\ge }{N-2}r^{-1}=C^*\ge r^{-1},
 $$
and 
$$m r^{1-N}\int_0^rs^{N-1}\overline{|\nabla u|^q}(s)ds\leq mC^q r^{1-N}\int_0^r|\ln s|^qs^{-q+N-1}ds \leq \tilde Cr^{1-q}.
$$
Hence
\bel{S14-1}
0\leq -\bar u_r(r)\leq C^*\ge r^{-1}+\tilde Cr^{1-q},
\ee
which finally implies that 
\bel{S15}
\bar u(r)\leq \bar u(r_\ge)+C^*\ge \ln \frac{r_\ge}{r}+\frac{\tilde C}{2-q}r_\ge^{2-q},
\ee
and
\bel{S16}
0\leq  e^{u(x)}\leq \left(e^{K'_5+\bar u(r\ge)+\frac{\tilde C}{2-q}r_\ge^{2-q}}r_\ge^{C^*\ge}\right) r^{-C^*\ge}.
\ee
Since $\ge$ is arbitrary, this implies that $e^{u}\in L^p(B_{r_0})$ for any $p>\frac N2$. Hence $\BBG_{B_{r_0}}[e^u]\in L^\infty(B_{r_0})$. Using $(\ref{Q12*})$ we conclude that $u$ remains bounded and then, by standard regularity estimates, we obtain that it coincides with a smooth solution in $B_{r_0}$.\smallskip

\nind \nind{\it Proof of \rth{supth6-bis}: Assertion 2}. Let us assume that $u$ satisfies $(\ref{Na-17})$. By assertion 1 this cannot hold if $q>\frac {N}{N-1}$; if $q=\frac{N}{N-1}$ $u$  is bounded from below. 
However by \rth{supth2}, $|\nabla u|^q$ would be locally integrable in $B_{r_0}$. Therefore $(\ref{Na-17})$ can hold only when $1<q<\frac{N}{N-1}$, and then $(\ref{Na-10})$ holds by \rth{supth1-1}.
Furthermore
\bel{S17}-(r^{N-1}\bar u_r)_r\geq mr^{N-1}\overline{|\nabla u|^q}(r).
\ee
Because of $(\ref{Na-17})$, for any $k>0$ there exists $r_k\in (0,r_0)$ such that $0<r\leq r_k$ implies  $\dsps\int_r^{r_0}\overline{|\nabla u|^q}(s)s^{N-1}ds>k$. By integrating $(\ref{S17})$ we obtain that
$$
r^{N-2}\bar u(r)\leq r^{N-2}\bar u(r_k)-\frac{mk}{N-2}\left(1-r^{N-2}r_k^{2-N}\right)-\frac{r_0^{N-1}\bar u_r(r_0)}{N-2}\left(1-r^{N-2}r_k^{2-N}\right)\quad\text{for }0<r\leq r_k.
$$
Hence $\dsps \limsup_{r\to 0}r^{N-2}\bar u(r)\leq -\frac{mk}{N-2}$. Since $k$ is arbitrary we obtain as a first result  that $r^{N-2}\bar u(r)\to -\infty$ when $r\to 0$, and this implies $(\ref{Na-17*})$.\\
If we assume that $u$ is bounded from above, then there exists $\gm>0$ such that the function $U_\gm=-u+\gm$ is bounded from below by $1$ and satisfies
 \bel{S17-1}
-\Gd U_\gm+m|\nabla U_\gm|^q+e^{\gm}e^{-U_\gm}=0.
\ee
Since $|\nabla U_\gm(x)|^{q-1}\leq c|x|^{-1}$ from $(\ref{Na-13})$, and $\frac{e^{\gm}e^{-U_\gm}}{U_\gm}\leq e^{\gm-1}$, Harnack inequality holds on every sphere of radius $r\leq \frac{r_0}{2}$  in the sense that
 \bel{S17-2}
\max_{|x|=r}U_\gm(x)\leq C\min_{|x|=r}U_\gm(x),
\ee
for some $C>0$. 
Combined with $(\ref{Na-17*})$ this inequality implies that for any sequence $\{r_n\}$ converging to $0$ there holds 
  \bel{S17-3}
\lim_{r_n\to 0}r_n^{N-2}U_\gm(r_n,.)\to\infty\quad \text{uniformly on }\, S^{N-1}.
\ee
By \rth{supth1} $\bar u(r)$ is a monotone function of $r$ and it tends to $-\infty$ when $r\to 0$. By Harnack inequality we  deduce that $U_\gm(x)\to \infty$ when $x\to 0$ and thus $\dsps \lim_{x\to 0}e^{-U_\gm(x)}= 0$. For $\ge>0$, there exists $r_\ge\in (0,r_0]$ such that $e^{-U_\gm(x)}\leq \ge$ in $B_{r_\ge}\setminus\{0\}$. This implies that $U_\gm$ satisfies 
 \bel{S17-1*}
-\Gd U_\gm+m|\nabla U_\gm|^q+\ge\geq 0\quad\text{in }\,B_{r_\ge}\setminus\{0\}.
\ee
For $m'>m$ and $\tilde\gg>0$ let $V_{\gg,m'}$ be the solution of 
  \bel{S17-4}
-\Gd V_{\tilde\gg,m'}+m'|\nabla V_{\tilde\gg,m'}|^q=c_N\tilde\gg\gd_0\quad\text{in }\CD'(\BBR^N),
\ee
which tends to $0$ when $|x|\to \infty$. This equation reduces to the separable first order one with the unknown $\frac{dV_{\tilde\gg,m'}}{dr}$
  \bel{S17-5}\left\{\BA{lll}\dsps
-\frac{d^2V_{\tilde\gg,m'}}{dr^2}-\frac{N-1}{r}\frac{dV_{\tilde\gg,m'}}{dr}+m'\left|\frac{dV_{\tilde\gg,m'}}{dr}\right|^q=0\quad\text{in }(0,\infty)\\[2mm]\phantom{-}
\dsps \lim_{r\to 0}r^{N-1}\frac{dV_{\tilde\gg,m'}(r)}{dr}=(1-N)\tilde\gg.
\EA\right.\ee
The solution is explicit, namely
   \bel{S17-6}\dsps
\frac{dV_{\tilde\gg,m'}}{dr}(r)=-\frac{r^{1-N}}{\left(((N-1)\tilde\gg)^{1-q}-\dsps\frac{m'}{\gk}r^{N-q(N-1)}\right)^{\frac 1{q-1}}},
\ee
(note that $\gk=q(N-1)-N<0$) and
   \bel{S17-7}\dsps
V_{\tilde\gg,m'}(r)=\int_r^\infty\frac{s^{1-N}}{\left(((N-1)\tilde\gg)^{1-q}-\dsps\frac{m'}{\gk}s^{N-q(N-1)}\right)^{\frac 1{q-1}}}ds,
\ee
which implies that
   \bel{S17-8}\dsps
V_{\tilde\gg,m'}(r)=\tilde\gg\frac{N-1}{N-2} r^{2-N}(1+o(1))\quad \text{as }\,r\to 0.
\ee
For any $m'>m$ we can choose $\ge$ and $r_\ge$ small enough such that for any $\tilde \gg\geq 1$,
$$(m'-m)\left|\frac{dV_{\tilde\gg,m'}}{dr}\right|^q\geq\ge\quad\text{in }B_{r_\ge}\setminus\{0\}.
$$
For such a choice, we have
  \bel{S17-9}
-\Gd V_{\gg,m'}+m|\nabla V_{\gg,m'}|^q+\ge\leq 0\quad\text{in }B_{r_\ge}\setminus\{0\}.
\ee
Next we combine $(\ref{S17-3})$, $(\ref{S17-8})$ and the maximum principle applied in $B_{r_\ge}\setminus B_{r_n}$ to $(\ref{S17-1*})$ and $(\ref{S17-9})$  to conclude that 
  \bel{S17-9bis}
U_\gm(x)\geq V_{\gg,m'}(|x|)-\max_{\prt B_{r_\ge}}U_\gm\quad\text{for all  }x\in B_{r_\ge}\setminus B_{r_n}.
\ee
Letting $r_n\to 0$ implies that the previous inequality is valid in $B_{r_\ge}\setminus \{0\}$, independently of $\gg$. Letting $\gg\to\infty$ we obtain that
$$U_\gm(x)\geq \int_{|x|}^\infty\left(-\frac{\gk)}{m'}\right)^{\frac{1}{q-1}}s^{-\frac{1}{q-1}}ds-\max_{\prt B_{r_\ge}}U_\gm=\frac{q-1}{2-q}\left(-\frac{\gk}{m'}\right)^{\frac{1}{q-1}}|x|^{\frac{q-2}{q-1}}-\max_{\prt B_{r_\ge}}U_\gm,
$$
in $B_{r_\ge}\setminus B_{r_n}$. If $m'\downarrow m$ and $r_n\to 0$
we obtain
$$\liminf_{x\to 0}|x|^{\gb}U_\gb(x)\geq \frac{q-1}{2-q}\left(-\frac{\gk}{m}\right)^{\frac{1}{q-1}},
$$
which implies $(\ref{Na-17**})$.
\qeda

\subsubsection{ Proof of \rth{supth6-ter}}
(I) {\it Step 1: case $u$ is decreasing near $0$}. Then $u$ is bounded from below. Up to changing $u$ into $u+\gm$ and modifying $m$ accordingly we can assume that $u(r)\geq 0$ when $r\to 0$. We perform the change of variable $(\ref{L4})$, with 
$$t=\ln r\,,\;x(t)=r^2e^{u(r)}\, \text{ and }\; \Gf(t)=-ru_r(r),$$ 
then 
                  \bel{S18}\left\{ \BA{lll}
x_t=x(2-\Phi)\\
\Phi_t= x-(N-2)\Phi+me^{(2-q)t}|\Phi|^q,
\EA\right.
   \ee
and $x$ is positive and bounded by $(\ref{Na-10})$ and $\Phi$ is positive. We aim to prove that $\Phi$ is also bounded when $t\to-\infty$. We define
\bel{S18*}\BA{lll}\dsps\CW_m(t)=:x\left(2(N-2)-\frac{(2-\Phi)^2}{2}-\frac{x}{2}-m2^qe^{(2-q)t}\right)-2(N-2)^2\\[2mm]
\phantom{\dsps\CW_m(t)}\dsps=-x\left[\frac{(2-\Phi)^2}{2}+m2^qe^{(2-q)t}\right]-\frac{(x-2(N-2))^2}{2}<0.
   \EA\ee
Then
   $$\BA{lll}\dsps\CW_{m\,t}=x\left[(2-\Phi)^2\left(N-3+\frac\Phi 2\right)-me^{(2-q)t}\left(2^q(2-q)+(2-\Phi)(2^q-\Phi^q\right)\right]\\[2mm]
   \phantom{\dsps\CW_{m\,t}}\dsps=x\left[(2-\Phi)^2\left(N-3+\frac\Phi 2\right)\right]+H,
   \EA$$
   where
   $$\frac H{xme^{(2-q)t}}=-2^q(2-q)-(2-\Phi)(2^q-\Phi^q)<0.
   $$
   We give below some properties of the function $\CW_{m}$. Since $x$ is bounded $\CW_m$ is bounded if and only if $x(2-\Phi)^2$ is bounded too. \smallskip

  \nind 
(1) First case, we assume that $\CW_m$ is bounded. We claim that $\Phi$ is bounded too. Assuming the contrary, we encounter two possibilities: \smallskip

\nind (1)-a The function $\Phi$ is unbounded and not monotone. Then there exists $\{t_n\}\to-\infty$ such that $\Phi(t_n)\to\infty$, $\Phi_t(t_n)=0$ and $\Phi_{tt}(t_n)\leq 0$. Hence\smallskip

   $x(t_n)-(N-2)\Phi(t_n)+me^{(2-q)t_n}|\Phi(t_n)|^q=0$ and 
   $x(t_n)(2-\Phi(t_n))+(2-q)me^{(2-q)t_n}|\Phi(t_n)|^q\leq 0,$\smallskip
   
   \nind which implies
   $$x(t_n)(2-\Phi(t_n))+(2-q)\left((N-2)\Phi(t_n)-x(t_n)\right)\leq 0,$$
   and thus
   $$qx(t_n)\leq\Phi(t_n)\left(x(t_n)-(2-q)(N-2)\right).
   $$
This in turn yields $x(t_n)>(2-q)(N-2)$ and therefore $x(t_n)(2-\Phi(t_n))^2\to\infty$ since $\Phi(t_n)\to\infty$, a fact which  is impossible by the assumption on $\CW_m$. \smallskip

\nind (1)-b The function $\Phi$ is unbounded and monotone  and thus it tends to 
   $\infty$. Thus $\frac {x_t(t)}{x(t)}\leq -A$ for $t\leq t_A<0$,$\!\!\!\!\!\!\!\phantom{\frac A{\int_a^b}}$ and clearly $x(t)\to\infty$ as $t\to-\infty$, which is also a contradiction. \\
  \nind {\it In conclusion $\CW_m$ is bounded if and only if $\Phi$ is bounded too.}\smallskip
   
  \nind (2) Second case, we  assume that $\CW_m$ is unbounded, but not monotone. Then there exists a sequence $\{\gt_n\}$ tending to $-\infty$ such that $\CW_{m\,t}(\gt_n)= 0$, $\CW_m(\gt_n)\to\-\infty$ and thus $x(\gt_n)(2-\Phi(\gt_n))^2\to\infty$. Set $\Phi(t_n)=2-h(\gt_n)>0$, 
then $x(\gt_n)h^2(\gt_n)\to\infty$. If we replace $\Phi(\gt_n)$ in the expression of $\CW_{m\,t}(\gt_n)= 0$ by $2-h(\gt_n)$, we obtain
\bel{S19*}h^2(\gt_n)\left(N-2-\frac{h(\gt_n)}2\right)=me^{(2-q)\gt_n}\left(2^q(2-q)+h(\gt_n)(2^q-(2-h(\gt_n))^q)\right).
   \ee
   Hence 
   $$h(\gt_n)\left(h(\gt_n)\left(N-2-\frac{h(\gt_n)}2\right)-me^{(2-q)\gt_n}\left(2^q-(2-h(\gt_n))^q)\right)\right)=me^{(2-q)\gt_n}2^q(2-q).
   $$
We have $N-2-\frac{h(\gt_n)}2=N-3-\frac{\Phi(t_n)}{2}>0$ and $h(\gt_n)(2^q-(2-h(\gt_n))^q)=(2-\Phi(t_n))(2^q-\phi(t_n)^q)>0$, we encounter two possibilities:\\
  - either $h(\gt_n)\geq 0$, then $N-2\geq\frac {h(\gt_n)}2$, then $0\leq h^2(\gt_n)\leq 4(N-2)^2$ which contradicts the fact that
   $x(\gt_n)h^2(\gt_n)\to\infty$. \\
- or $h(\gt_n)<0$ and from $(\ref{S19*})$ we have
   $$\BA{lll}\dsps
   \frac{|h(\gt_n)|^3}{2}\leq h^2(\gt_n)\left(N-2+\frac{|h(\gt_n)|}2\right)\\
   \phantom{\dsps\frac{|h(\gt_n)|^3}{2}}\dsps\leq 
   me^{(2-q)\gt_n}\left[2^q(2-q)+|h(\gt_n)|((2+|h(\gt_n)|)^q-2^q)\right]
   \\[0mm]
   \phantom{\dsps\frac{|h(\gt_n)|^3}{2}}\dsps\leq me^{(2-q)\gt_n}(c(q)+|h(\gt_n)|)^{q+1}).
   \EA$$
Therefore $|h(\gt_n)|\leq ce^{(2-q)\gt_n}$ since $q<2$ which again contradicts  the fact that $x(\gt_n)h^2(\gt_n)\to\infty$. \smallskip

\nind {\it Hence the second case never holds.}

\nind (3) Third case, we  assume that $\CW_m$ is unbounded, and monotone, then
 $\CW_m(t)$ tends to $-\infty$. The expression of $\CW_m$ shows that $\Phi(t)\to\infty$ when $t\to\infty$, hence $\frac {x_t}x\to-\infty$ as $t\to-\infty$, which is not compatible with the boundednes of $x$. \smallskip

\nind {\it Hence the third case never occurs.}\smallskip

\nind {\it If the function $u$ is decreasing, the function $\CW_m$ is bounded, and then $x$ and $\Phi$ are bounded too.}\medskip

The system $(\ref{S18})$ is an exponential perturbation of the simpler system associated to the equation $-\Gd u=e^u$,
                  \bel{S19} \left\{\BA{lll}
x_t=x(2-\Phi)\\
\Phi_t= x-(N-2)\Phi,
\EA\right.
\ee
which admits two equilibria in $\BBR^2$, $(0,0)$ and $(2(N-2),2)$. The point $(0,0)$ is a saddle point with characteristic values $2$ and $2-N$, and the point $(2(N-2),2)$ is a sink.
The function $\CW_0$ defined in $(\ref{S18*})$ is nonpositive and increasing. Therefore the system $(\ref{S19})$ admits no periodic solution in the quadrant $\{x>0,\Phi>0\}$. Since 
 $x$ and $\Phi$ are bounded we can apply \cite[Theorem 4.1]{LoRy} (see also \cite[Proposition 13]{BV-V25}) the limit set at $-\infty$ of any 
bounded solution of $(\ref{S18})$ is invariant by the flow of the system $(\ref{S19})$ (always in the same quadrant). The only invariant sets for this flow are $(0,0)$ and $(2(N-2),2)$. therefore, either
the solution $(x(t),\phi(t))$ converges to $(2(N-2),2)$ when $t\to-\infty$, which is exactly $(\ref{Na-17-2})$ or it converges to $(0,0)$, and this is more delicate.\\
In such a case we are led to introduce the system of order 3 where $\Gth(t)=e^{(2-q)t}$,
                  \bel{S20} \left\{\BA{lll}
x_t=x(2-\Phi)\\
\Phi_t= x-(N-2)\Phi+m\Phi^q\Gth\\
\Gth_t=(2-q)\Gth.
\EA\right.
\ee
Notice that $\Gth(t)\to 0$ when $t\to-\infty$. Since $(2(N-2),2)$ is a sink, standard perturbation theory asserts that 
there exists a neighbourhood $\CV_1$ of this point such that all the trajectories of $(\ref{S18})$ issued from this neighbourhood converge to $(2(N-2),2)$ when $t\to-\infty$. This means that relation $(\ref{Na-17-2})$ in \rth{supth6-ter} holds. \\The analysis of $(0,0)$ is more delicate since the system $(\ref{S20})$ admits the axis $x=\Phi=0$ as a trajectory. The linearisation of $(\ref{S20})$ at $(0,0,0)$ is 
                  \bel{S21} \left\{\BA{lll}
x_t=x(2-\Phi)\\
\Phi_t= x-(N-2)\Phi\\
\Gth_t=(2-q)\Gth.
\EA\right.
\ee
The eigenvalues of the system are $\gl_1=2>0$, $\gl_2=2-N<0$ and $\gl_3=2-q>0$. Then there exists a 2-dimensional stable manifold
$\CM$ of $\BBR^3$ of solutions converging to $(0,0,0)$ when $t\to-\infty$. The tangent plane to $\CM$ at $(0,0,0)$ is spanned by the two eigenvectors $\gw_1=(N,1,0)$ and $\gw_3=(0,0,1)$. Now the system $(\ref{S20})$ reduced to the plane 
$x=0$ is 
\bel{S22} \left\{\BA{lll}
\Phi_t= x-(N-2)\Phi+m\phi^q\Gth\\
\Gth_t=(2-q)\Gth.
\EA\right.
\ee
The eigenvalues of its linearisation at $(0,0)$ are $\gl_2=2-N<0$ and $\gl_3=2-q>0$, hence there exists only  one stable curve of solutions converging to $(0,0)$ when $t\to-\infty$, and it is precisely $\Phi=0$ and $\Gth(t)=e^{(2-q)t}$. It means that $\CM$ is not contained into the plane $x=0$, hence the trajectories contained in $\CM\cap\{x>0\}$ are admissible in the sense that they correspond to solutions of $(\ref{S18})$. Since $0<\gl_3<\gl_1$, there exist infinitely many trajectories of $(\ref{S20})$ contained in $\CM$ admitting $\gw_3$ for tangent vector at $(0,0,0)$ (all except the trajectory in the plane $x=0$). These trajectories which are tangent to the axis $x=\Phi=0$ at $0$ satisfies 
$$x(t)+\Phi(t)=r^2e^{u(r)}+r|u_r(r)|=o(e^{(2-q)t}).$$ 
Hence $u(r)=o(\ln\frac 1r)$ and $|u_r(r)|=o(r^{1-q})$, which implies that 
$u_r$ is integrable, hence there exists $u_0$ such that $u(r)\to u_0$ as $r\to 0$, and $x(t)=e^{u_0}e^{2t}(1+o(1))$. Moreover, since we have assumed that $u$ is decreasing, we have $\Phi>0$. Finally,  since $\Phi(t)\to 0$, there holds 
$$\Phi_t<x-\frac {N-2}{2}\Phi<x=e^{2t}e^u<e^{2t}e^{u_0}=\frac 12(e^{2t}e^{u_0})_t.
$$
This implies that $t\mapsto \Phi(t)-\frac12e^{2t}e^{u_0}$ is decreasing. Since it tends to $0$, we have 
$\Phi(t)\leq \frac12e^{2t}e^{u_0}$ and therefore $|u_r|\leq \frac 12re^{u_0}$. This implies that $u_r(r)\to 0$ when $r\to 0$, which means that $u$ is a regular solution. As a consequence $u$ satisfies either
$(\ref{Na-17-2})$ or $(\ref{Na-17-3})$. Note that the fact that $u_r(r)=o(r^{1-q})$ 
can also be obtained by using the non-radial result of \rth{supth5} proved by the Bernstein technique.
\medskip

\nind (II) {\it Step 2: case $u$ is increasing near $0$}. Then either $u(r)$ admits a finite limit $u_0$ or 
it tends to $-\infty$. Moreover $u_r$ is monotone near $0$. Indeed at a point $r^*$ near $0$ where $u_{rr}(r^*)=0$, we have
$$u_{rrr}(r^*)=\left(\frac{N-1}{r^{*2}}-e^{u(r^*)}\right)u_r(r^*)>0,
$$
since either $e^{u(r)}\to e^{u_0}$ or $e^{u(r)}\to 0$. Then if $u$ is bounded,\smallskip

\nind - Either $u_r\to 0$ when $r\to 0$, and then $u$ is a regular, but in that case we have $-Nu_{rr}((r)\to e^{u_0}$, which implies that $0$ is a local maximum of $u$. By Taylor-Young expansion $u_r(r)<0$ for $r$ small enough which is a  contradiction.\

\nind - Or $u_r(r)\to c_0>0$, and then $u_{rr}(r)=\frac{1-N}{r}c_0(1+o(1))$ which implies that $u_{rr}$ is not integrable at 
$0$ and $u_r$ has no limit at $0$ which is a  contradiction. \smallskip

\nind It follows that $u_r(r)\to\infty$, then $u$ is negative near $0$. Setting $U=-u$, then $U>0$ and 
$$- U_{rr}-\frac{N-1}{r}U_r+e^{-U}+m|U_r|^q=0.
$$
Setting $W=-r^{N-1}U_r>0$, then $W_rW^{-q}+mr^{(N-1)(1-q)}\leq 0$. Hence, if $q\neq \frac{N}{N-1}$, the function $r\to\Psi(r)=\frac{W^{1-q}}{q-1}+\frac{m}{(N-1)q-N}r^{N-(N-1)q}$ is increasing. 
If $q>\frac{N}{N-1}$ this is impossible. If $q=\frac{N}{N-1}$  the function $r\to (N-1)W^{\frac 1{1-N}}-m\ln r$ is increasing, which also is impossible.
Then necessarily $1<q<\frac N{N-1}$ and since $e^{-U}\leq 1$, we have that $e^{-U}=o(|U_r|^q)$, thus $U$ satisfies the following viscous Hamilton-Jacobi equation,
$$- U_{rr}-\frac{N-1}{r}U_r+\tilde m(r)|U_r|^q=0,
$$
where $\tilde m(r)=(m+e^{-U}|U_r|^{-q})=m(1+o(1))$. The function $W$ defined above satisfies
$$\frac{d}{dr}\left(\frac{W^{1-q}(r)}{1-q}+\int_0^r\tilde m(s)s^{(N-1)(1-q)}ds\right)=0.
$$
This identity implies that 
\bel{S23} 
\frac{W^{1-q}(r)}{1-q}+\int_0^r\tilde m(s)s^{(N-1)(1-q)}ds=\ell
\ee
for some real number $\ell $. Since
$$\int_0^r\tilde m(s)s^{(N-1)(1-q)}ds=\frac{m}{N-(N-1)q}r^{N-(N-1)q}(1+o(1))\quad\text{when }r\to 0,
$$
we have the following identity if $\ell=0$:
\bel{S24} 
W^{1-q}(r)=-\frac{m}{\gk}r^{N-q(N-1)}(1+o(1)),
\ee
which implies
$
U_r(r)=-\left(-\frac{\gk}{m}\right)^{\frac1{q-1}}r^{-\frac{1}{q-1}}(1+o(1))
$
and
\bel{S25} 
u(r)=-\frac{2-q}{q-1}\left(-\frac{\gk}{m}\right)^{\frac1{q-1}}r^{-\frac{2-q}{q-1}}(1+o(1))\quad\text{when }r\to 0, 
\ee
which is $(\ref{Na-17-4})$. Next, if $\ell\neq 0$, then clearly $\ell<0$,
$$W(r)=\left(\frac 1{q-1}\right)^{\frac 1{q-1}}|\ell|^{-\frac{1}{q-1}}+o(1),
$$
and 
\bel{S26} 
u(r)=-\left(\frac{1}{(q-1)|\ell|}\right)^{\frac1{q-1}}r^{2-N}+o(1)\quad\text{when }r\to 0, 
\ee
which is $(\ref{Na-17-5})$ with $\gg=-\left(\frac{1}{(q-1)|\ell|}\right)^{\frac1{q-1}}$.\qeda

\subsubsection{ Proof of \rth{supth7} and \rth{supth7bis}}

\medskip

The following counterpart of \rlemma{isotropy} which is  a variant  of \cite[Proposition 5.6]{BV-V24} holds.

\blemma{isotropyinf} Let $N\geq 2$, $q>2$ and $u$ be a solution of $(\ref{Na-1})$ in $B_{r_0}^c$ such that $|x|^2e^u$ is bounded.
 Then there exists $C_1>0$ depending on $N, p, m$ and $\norm{|x|^2e^{u}}_{L^\infty}$ such that 
\bel{S27}
\norm{u(r,.)-\bar u(r)}_{L^\infty(S^{N-1})}\leq C_1\quad\text{for all }r\geq 2r_0.
\ee
\es
The variation lies in the sign of $m$, but this plays no role in the application of the representation formula\cite[formula (5.20)]{BV-V24}.  \medskip

The results of \rlemma{Reg} are also valid if we replace $t\geq t_0$ by $t\leq t_0$, provided we assume $q>2$ (hence $\gb<0$) with a function $\eta\in C^2((-\infty,t_0])$ which satisfies  $(\ref{Q6})$. Then 
$(\ref{Q7})$ is valid for $(t,\gs)\in (-\infty,t_0]\ti S^{N-1}$.\medskip

\nind{\it Proof of \rth{supth7}}. Under the assumption that $|x|^2e^{u(x)}\leq M$ the function $u$ is upper bounded. Combining $(\ref{Na-11***})$ in   \rth{supth3} with  (\ref{Na-12-bis}) in \rcor{supth4bis} we conclude that $u$ satisfies 
\bel{S27*}
-C|x|^{-\gb}\leq u(x)\leq C|x|^{-\gb}\quad\text{for }|x|\geq 2r_0.
\ee

\nind 1- If there holds
\bel{S28}
0<\liminf_{r\to\infty}|x|^2e^{u(x)}\leq\limsup_{r\to\infty}|x|^2e^{u(x)}<\infty,
\ee
we adapt the proof of \rth{supth6} to the exterior domain and set
\bel{QS2}r=e^{t}\Longleftrightarrow t=\ln r\quad\text{and }\,v(t,\gs)=u(r,\gs)-2\ln\frac 1 r\Longrightarrow v(t,\gs)=u(r,\gs)+2t,\ee
and there holds in $(t_0,\infty)\ti S^{N-1}$ (with $t_0=\ln r_0$)
\bel{QS3}v_{tt}+(N-2)v_t+\Gd'v-2(N-2)+e^{v}+me^{(2-q)t}\left(\left(v_t-2\right)^2+|\nabla 'v|^2\right)^{\frac q2}=0.
\ee
As $u$ is a solution of $(\ref{Na-1})$ in $B^c_{r_0}$ such that $|x|^2e^u\leq\tilde K_1$, the function $v$ satisfies 
\bel{Q4-bis} v(t,\gs)\leq  K_1.
\ee
Moreover, since the function $u(x)-2\ln\frac1{|x|}$ is bounded on $B_{2r_0}^c$, then the functions $v$ and consequently $v_t$, $\nabla'v$, $v_{tt}$, $\nabla'v_t$ and $\nabla'^2v$ are also bounded by the same argument as in \rlemma{Reg}. The dynamical system approach and the utilisation of 
Huang-Takac and Simon's results \cite{Sim}, \cite{HuTak} as in the proof of \rth{supth6} applies. By \rlemma{isotropyinf} estimate $(\ref{S27})$ holds, then we obtain $(\ref{Na-21*})$.\smallskip

\nind 2- Next we assume that $\liminf_{r\to\infty}|x|^2e^{u(x)}=0$. Then for any $k>0$ there exists $x_k$ with $r_k:=|x_k|>r_0$ such that $u(x_k)\leq -k-2\ln r_k$. Using \rlemma{isotropyinf} we deduce that 
$u(x) \leq -k-2\ln r_k+C_1$ for any $x$ such that $|x|=r_k$.
We set $U=-u$, then 
\bel{S28*}-\Gd U+m|\nabla U|^q+e^{-U}=0,
\ee
and $U(x)\geq k+2\ln r_k-C_1$ for all $|x|= r_k$, and $r_k\to\infty$ when $k\to\infty$. Furthermore $\bar U(r)$ is increasing by \rth{supth1} in $B_{r_1}^c$ for some $r_1>r_0$. Let $\tilde m>m$. For $a>0$ we introduce the problem
\bel{S29}\left\{\BA{lll}
-W_r+\tilde mr^{(N-1)(1-q))}W^q=0\quad\text{in }(r_0,\infty)\\
\phantom{+mr^{(N-1)(1-q))}W^q}
W(r_0)=a.
\EA\right.\ee
The solution denoted by  $W_{a,\tilde m}$ is increasing on the maximal interval $[r_0,r_a)$, with $r_a\leq\infty$, with limit $L_{a,\tilde m}\in (a,\infty]$ at $r=r_a$, and we have
$$\frac d{dr}\left[ W_{a,\tilde m}^{1-q}(r)-\frac{\tilde m}{\gk}r^{N-q(N-1)}\right]=0.
$$
If 
\bel{S30}
a<\left(\frac{\gk}{\tilde m}\right)^{\frac 1{q-1}}r_0^{\frac{q(N-1)-N}{q-1}}:=a^*
\ee
then $r_a=\infty $ and $L_{a,\tilde m}<\infty$, if $a=a^*$ (resp. $a>a^*$) then $r_a=\infty$ and $L_{a^*,\tilde m}=\infty$ (resp. $r_a<\infty$ and $L_{a,\tilde m}=\infty$. Note also that 
\bel{S31}
W_{a^*,\tilde m}(r)=\left(\frac\gk{\tilde m} \right)^{\frac 1{q-1}}r^{\frac{q(N-1)-N}{q-1}}.
\ee
With this expression we can express the solution of $V=V_{A,a^*,\tilde m}$ of 
\bel{S33}\left\{\BA{lll}
-V_{rr}-\frac{N-1}{r}V_r+\tilde m|V_r|^q=0\quad\text{in }(r_0,\infty)\\[1mm]
\phantom{+\tilde mr^{(N-1)(1-q))}W^q}
\!V(r_0)=A,\, V_r(r_{0})=a^*.
\EA\right.\ee
Using
$V_r(r)=r^{1-N}W_{a,\tilde m}(r)$, we obtain
\bel{S34}\BA{lll}\dsps
V_{A,a^*,\tilde m}(r)=A-\frac{1}{\gb}\left(\frac\gk{\tilde m} \right)^{\frac 1{q-1}}(r^{-\gb}-r_0^{-\gb}),
\EA\ee
in which expression we recall that $-\gb=\frac{q-2}{q-1}>0$. The last thing to insure is that $V_{A,a^*,\tilde m}$ is a subsolution in the domain $B_{\gr}\setminus B_{r_0}$ we aim to compare 
$U$ and $V_{A,a^*,\tilde m}$. We write the two equations that these functions satisfy under the form
$$\BA{lll}-\Gd U+m|\nabla U|^q=-e^{-U}\\[1mm]
-\Gd V_{A,a^*,\tilde m}+m|\nabla V_{A,a^*,\tilde m}(x)|^q=(m-\tilde m)|\nabla V_{A,a^*,\tilde m}|^q.
\EA$$
Since $|x|^2e^{-U(x)}\leq C^*$ and $|\nabla V_{A,a^*,\tilde m}(x)|^q=\left(\frac{\gk}{m}\right)^{\frac q{q-1}}|x|^{-\frac{q}{q-1}}$, we have 
$$(m-\tilde m)|\nabla V_{A,a^*,\tilde m}(x)|^q+e^{-U(x)}\leq (m-\tilde m)|x|^{-\frac{q}{q-1}}+C^*|x|^{-2}=|x|^{-2}\left(C^*-(\tilde m-m)|x|^{-\gb}\right)
$$
If $\tilde m>m$ is given, we can change $r_0$ (which does not impact the value of $|\nabla V_{A,a^*,\tilde m}(x)|$) is order to have 
$C^*<(\tilde m-m)r_0^{-\gb}
$
and consequently 
$$-\Gd V_{A,a^*,\tilde m}+m|\nabla V_{A,a^*,\tilde m}(x)|^q\leq -\Gd U+m|\nabla U|^q\quad\text{in }B_{r_0}^c.
$$
Next we take $A>0$ in order $A< U(r_0)$. 
For $k>1$ large enough such that $r_k>r_0$ and for $|x|=r_k$, we have $U(x)\geq \ln k+2\ln r_k-C_1$. Then:\\
- either $V_{A,a^*,\tilde m}(r_k)\leq \ln k+2\ln r_k-C_1$ for all $k$ large enough, in which case we deduce by the comparison principle that 
\bel{S34-1}\BA{lll}\dsps
V_{A,a^*,\tilde m}(|x|)\leq U(x)\quad\text{for all }x\in B_{r_0}^c,
\EA\ee
- or there exists a subsequence $\{k_j\}$ tending to infinity such that  
\bel{S34-2}\BA{lll}\dsps
V_{A,a^*,\tilde m}(|x|)\leq U(x)\quad\text{for all }\,x\in B_{r_{k_j}}\setminus B_{r_0}^c\\
V_{A,a^*,\tilde m}(|x_j|)> U(x_j)\quad\text{for some }\,x_j\in B_{r_{k_{j+1}}}\setminus B_{r_{k_j}}.
\EA\ee
Since $\{k_j\}\to\infty$ we deduce also that $(\ref{S34-1})$ holds. As a consequence
\bel{S34-3}\BA{lll}\dsps
-\frac{1}{\gb}\left(\frac\gk {\tilde m}\right)^{\frac 1{q-1}}=\lim_{|x|\to\infty}|x|^{\gb}V_{A,a^*,\tilde m}(|x|)\leq \liminf_{|x|\to\infty}|x|^{\gb}U(x)
\EA\ee
Since $\tilde m>0$ is arbitrary, we let it converge to $m$. Combining this resulting inequality with $(\ref{Na-11**})$ we obtain $(\ref{Na-22})$.
\qeda
\medskip

\nind{\it Proof of \rth{supth7bis}}. From inequality $(\ref{Na-10})$ in  \rth{supth2} the function $r^2e^{u(r)}$ is bounded. Then \rth{supth7} implies \rth{supth2}. However we give below a  proof
specific to the radial case which puts into light some beautiful constructions originated in the theory of finite dimensional dynamical systems.\smallskip

We set $U(r)=-u(r)$, $\xi(t)=r^{\gb}U(r)$, $\eta (t)=r^{\frac 1{q-1}}U_r(r)$ and $t=\ln r$. The functions $\xi$ and $\eta$ are positive since $U(r)\to\infty$ and $U_r(r)$ is positive. As $U$ is a radial solution of $(\ref{S28*})$ , $(\xi,\eta)$ satisfies 
\bel{T3}\left\{\BA{lll}
\xi_t=\gb\xi+\eta\\
\eta_t=-\gk\eta+m|\eta|^q+e^{\frac{q}{q-1}t-e^{-\gb t}\xi},
\EA\right.\ee
where we recall that $\gk=\frac{q(N-1)-N}{q-1}$ and $\gb=\frac{2-q}{q-1}<0$. We set 
$$G(t)=e^{-\gb t}\xi(t)-\frac{q}{q-1}t=U(r)-\frac{q}{q-1}t,$$
and $(\xi,\eta)$ is bounded at infinity from $(\ref{Na-11})$ and $(\ref{Na-11**})$. By $(\ref{Na-12})$, $U(r)\geq 2t-C$, hence
$G(t)\geq  \frac{q-2}{q-1} t-C$. Therefore $(\ref{T3})$ is an exponential perturbation of the system associated to the Hamilton-Jacobi equation $(\ref{Na-2})$,
\bel{T4}\left\{\BA{lll}
\xi_t=\gb\xi+\eta\\
\eta_t=-\gk\eta+m|\eta|^q,
\EA\right.\ee
which admits two equilibria $(\xi_m,\frac{q-2}{q-1}\xi_m)$ and $(0,0)$ with $\xi_m=\frac{q-1}{q-2}\left(\frac\gk m\right)^\frac{1}{q-1}$. The point $(\xi_m,\frac{q-1}{q-2}\xi_m)$ is a saddle point, 
 while the origin is a sink This system admits no closed orbit (since the axis $\{\eta=0\}$ contains trajectories). Therefore the two equilibria are the only connected and compact invariant sets for $(\ref{T4})$. By the results of \cite[Theorem 4.1]{LoRy}, we have that  $(\xi(t),\eta(t))$ converges, either to $(0,0)$ or $(\xi_m,\frac{q-2}{q-1}\xi_m)$ to when $t\to\infty$. If $(\xi(t),\eta(t))\to (\xi_m,\frac{q-2}{q-1}\xi_m)$ we obtain in particular $(\ref{Na-22})$ and $\dsps\lim_{r\to \infty}r^{\frac 1{q-1}}u_r(r)=\left(\tfrac\gk m\right)^\frac{1}{q-1}$. \\
Now, if $(\xi(t),\eta(t))\to (0,0)$ the problem is more delicate. We write $(\ref{T3})$ under the form
\bel{T5}\left\{\BA{lll}
\xi_t=\gb\xi+\eta\\
\eta_t=-\gk\eta+m\eta^q+c(t)e^{\gb t},
\EA\right.\ee
where $c(t)$ is a bounded positive function. Since $\eta(t)\to 0$, for any $\ge>0$ we have
$$\eta_t\leq (-\gk+\ge)\eta+c(t)e^{\gb t}
$$
for $t\geq t_\ge >t_0$. Hence
$$\eta(t)\leq e^{(-\gk+\ge)(t-t_\ge)}\eta(t_\ge)+e^{(-\gk+\ge)t}\int_{t_\ge}^te^{(\gk-\ge+\gb)s}c(s) ds.
$$
Notice that $\gk>-\gb$ since $q>1$. Hence the previous inequality implies $\eta(t)\leq Ce^{\gb t}$ with $C>0$ depending on $t_\ge$ and this finally yields
\bel{T6} 
0\leq -u_r(r):=U_r(r)\leq \frac Cr\,\text{ which implies }\;u(r)\geq -C\ln \frac{r}{r_0}+u(r_0).
\ee 
However this estimate is not enough to conclude. We use system $(\ref{S18})$ with $x(t)=r^2e^{u(r)}$ and $\Phi(t)=-ru_r$. Since $\Phi$ is bounded by 
$(\ref{T6})$ and $q>2$, this system is an exponential perturbation of  the system $(\ref{S19})$ which admits two equilibria in the domain $\{(x,\phi):x\geq 0,\phi\geq 0\}$; they are $(0,0)$ and 
$(2(N-2),2)$, the nature of which is studied in the proof of \rth{supth6-ter}: we recall that $(2(N-2),2)$ is a sink, while $(0,0)$ is a saddle point and there exists no periodic solution exists. \\
Therefore we can use  again the result of Logemann and Ryan \cite[Theorem 4.1]{LoRy} to conclude that any trajectory $(x(t),\Phi(t))$ converges either to $(0,0)$ or to $(2(N-2),2)$ when $t\to\infty$. This is equivalent to:\\
(i) either  $(r^2e^{u(r)},-ru_r(r))\to (2(N-2),2)$ when $r\to\infty$,\\
(ii) or $(r^2e^{u(r)},-ru_r(r))\to (0,0)$ when $r\to\infty$.\\
In case (i) this yields
$$u(r)=-2\ln r+\ln(2(N-2)+o(1)\,\text{ and }\;ru_r(r)=-2+o(1)\quad\text{as }r\to\infty.
$$
In case (ii) we deduce from $(\ref{S18})$ that 
$\frac{x_t}{x}:=(\ln x)_t=2+o(1)$, which implies $\ln(x(t))=2t+o(t)$ when $t\to\infty$. This is equivalent to 
$$u(r)+2\ln r=2\ln r+o(\ln r).
$$
Then $u(r)=o(\ln r)$ which is not compatible with $u(r)\leq -2\ln r+C$. Hence only case (i) occurs. This ends the proof.\qeda

\mysection{Existence of solutions}


\subsection{Proof of \rth{supth8}}

\subsubsection{Proof of assertion 1}
We use the variable $x(t)=r^2e^{u(r)}$, $\Phi(t)=-ru_r(r)$ and $t=\ln r$. Hence $(x,\Phi)$ satisfies the system $(\ref{S18})$ which is an exponential perturbation of the autonomous system 
$(\ref{S19})$. This last system admits $(2(N-2),2)$ among its equilibria. This point is a sink with real negative eigenvalues if $3\leq N\leq 10$ and complex eigenvalues with negative real part if 
$N\geq 11$. It is therefore classical (see e.g. \cite{Dieu}) that there exists $t_0>0$ and a neighbourhood $\CV$ of $(2(N-2),2)$ such that for any $(x_0,\Phi_0)\in \CV$ the positive trajectory of $(\ref{S18})$
starting from this point at $t=t_0$ converges to $(2(N-2),2)$ as $t\to\infty$. Equivalently the radial solution $u$ of $(\ref{Na-1})$ which satisfies $u_r(r_0)=\ln x_0-2\ln r_0$ and $u_r(r_0)=r_0^{-1}\Phi_0$ satisfies  $(\ref{Na-21})$ provided $r_0$ is large enough and $|x_0-2(N-2)|+|\Phi_0-2|\leq \gd$ for some $\gd>0$.\qeda

\subsubsection {Proof of assertion 2} 
We recall that the function $V_{A,a^*,\tilde m}$ defined in $(\ref{S34})$ satisfies $(\ref{S33})$ in $(r_0,\infty)$ and $V_{A,a^*,\tilde m}(r_0)=A$. We have 
\bel{Z1}
-\Gd V_{A,a^*,\tilde m}+m|\nabla V_{A,a^*,\tilde m}|^q+e^{-V_{A,a^*,\tilde m}}=(m-\tilde m)|\nabla V_{A,a^*,\tilde m}|^q+e^{-V_{A,a^*,\tilde m}}.
\ee
From the expression of $V_{A,a^*,\tilde m}$ we have
\bel{Z2}(m-\tilde m)|\nabla V_{A,a^*,\tilde m}|^q+e^{-V_{A,a^*,\tilde m}}=(m-\tilde m)\left(\frac{\gk}{\tilde m}\right)^{\frac q{q-1}}r^{-\frac q{q-1}}+e^{-A-\frac{1}{\gb}r_0^{-\gb}}e^{\frac{1}{\gb}r^{-\gb}}
\ee
We define $A^*=A^*(r_0)$ by
$$(\tilde m-m)\left(\frac{\gk}{\tilde m}\right)^{\frac q{q-1}}=e^{-A^*-\frac{1}{\gb}r_0^{-\gb}}:=K_{A^*}.
$$
Denote $M(r)=e^{\frac{1}{\gb}r^{-\gb}}-r^{-\frac q{q-1}}$. If $\gr=-\frac{1}{\gb}r^{-\gb}=\frac{q-1}{q-2}r^{\frac{q-2}{q-1}}$  we define $\tilde M$ by 
$$M(r)=\tilde M(\gr)=e^{-\gr}-\left(\frac{q-2}{q-1}\right)^{\frac{q}{q-2}}\gr^{\frac{q}{2-q}}.$$
There exists $\gr_1=\gr_1(q)$ such that for any $\gr\geq \gr_1$ there holds $\tilde M(\gr)\leq 0$, equivalently 
\bel{Z3}
M(r)\leq 0\quad\text{for all }\,r\geq r_1:=r_1(q)=\left(\frac{q-2}{q-1}\gr_1(q)\right)^{\frac{q-1}{q-2}}
\ee
Now we fix $r_0=r_1$ and for any $A\geq A^*(r_1)$ and $r\geq r_1$ we have
\bel{Z4}
(m-\tilde m)|\nabla V_{A,a^*,\tilde m}|^q+e^{-V_{A,a^*,\tilde m}}\leq 0.
\ee
Since $V_{A,a^*,\tilde m}$ and $V_{A,a^*, m}$ are respectively a subsolution and a supersolution of  $(\ref{Na-1})$, they are Lispchitz continuous 
and bounded in $B_{r_1}^c$ and satisfies $V_{A,a^*,\tilde m}\leq V_{A,a^*, m}$, for any $n>r_1$ there exists a viscosity solution $U_n\in L^\infty(\Gw_n)\cap W^{1,q}(\Gw_n)$ of 
$(\ref{Na-1})$ in $\Gw_n:=B_n\setminus \overline B_{r_1}$ radially symmetric as $V_{A,a^*,\tilde m}$ and $V_{A,a^*, m}$ are, and satisfying 
$$
V_{A,a^*,\tilde m}\leq U_n\leq V_{A,a^*, m}\quad\text{in }\Gw_n,
$$
(see e.g. \cite{Li}, \cite{DaP}). Furthermore since $U_n$ is radial, $r\mapsto U_{n\,r}(r)$ is uniformly continuous on $[r_1,n]$. It follows directly from the differential equation that $U_n$ belongs to $W^{1,\infty}(\Gw_n)\cap C^2(\Gw_n)$. We can let 
$n\to\infty$ and derive that there exist a subsequence, $U_{n_k}$ and a function $U$ in $L_{loc}^\infty(B^c_{r_1})\cap W^{1,q}_{loc}(B^c_{r_1})$, such that $U_{n_k}$ converges to 
$U$ locally in $B^c_{r_1}$ and weakly in $W^{1,q}_{loc}(B^c_{r_1})$. The function $U$ is a radial viscosity solution of $(\ref{Na-1})$ in $B^c_{r_1}$, it satisfies 
$$
V_{A,a^*,\tilde m}\leq U\leq V_{A,a^*, m}\quad\text{in }B^c_{r_1}.
$$
From the equation, the function $U$ belongs to $C^2(B^c_{r_1})$. The function $u=-U$ is a solution of $(\ref{Na-1})$ in $B_{r_1}^c$ with value $-A$ on $\prt B_{r_1}$ and by \rth{supth7bis} it satisfies 
$(\ref{Na-22})$. The gradient estimate $(\ref{Na-24-bis})$ is standard from $(\ref{Na-22})$.\qeda
\subsection{Proof of \rth{supth9}} 
The proof below is an adaptation of the construction presented in \cite[Theorem 15]{BV-V25}.
\subsubsection { Proof of Assertion 1} We use system $(\ref{S18})$, 
$$\left\{\BA{lll}x_t=x(2-\Phi)\\
\Phi_t=x-(N-2)\Phi+me^{(2-q)t}|\Phi|^q,
\EA\right.
$$
 with 
$$x(t)=r^2e^{u(r)}\,,\;\Phi(t)=-ru_r(r)\,,\;t=\ln r$$ 
and its autonomous extension $(\ref{S20})$ to $\BBR^3$
which is 
$$\left\{\BA{lll}
x_t=x(2-\Phi)\\
\Phi_t=x-(N-2)\Phi+m\Gth|\Phi|^q\\
\Gth_t=(2-q)\Gth.\EA\right.
$$
 The solution we look for satisfies $x(t)\to 2(N-2)$, $\Phi(t)\to 2$ and $\Gth(t)\to 0$ when $t\to-\infty$ with 
equilibria $P_0=(2(N-2),2,0)$ and $O=(0,0,0)$. We write $x=2(N-2)+\bar x$ and $\Phi=2+\bar\Phi$. The linearised system at $P_0$ is 
\bel{U1}\left\{\BA{lll}\dsps
\bar x_t=-2(N-2)\bar \Phi\\
\bar\Phi_t=\bar x-(N-2)\bar\Phi+2^qm\Gth\\
\Gth_t=(2-q)\Gth,
\EA\right.\ee
with characteristic polynomial 
$$P(\gl)=(2-q-\gl)(\gl^2+(N-2)\gl+2(N-2)).
$$
Its roots are $\gl_1=2-q>0$ and $\gl_2,\gl_3$ which could be real and distinct if $N>10$, real and equal if $N=10$ and complex non-real with real part $2-N$ if $3\leq N<10$.  An eigenvector associated to $\gl_1$ is $\dsps\gw_1=\left(2(N-2),q-2,f(q)\right)$ where $f(q)=q^2-(N+2)q+4(N-1)$. Standard verifications show that $f(q)\neq 0$ for $1<q<2$. hence there exists a unique trajectory 
$\CT_1$ of the system $(\ref{S20})$ converging to $P_0$ when $t\to-\infty$, admitting $\gw_1$ as tangent vector at $P_0$ and such that $\Gth(t)>0$ near $P_0$. Let $(x(t),\Gf(t),\Gth(t))$ be a solution of $(\ref{S20})$ with trajectory $\CT_1$. Then $\Gth_t=(2-q)\Gth$, thus $\Gth(t)=be^{(2-q)t}$ for some $b>0$. Therefore
$$\left\{\BA{lll}
x_t=x(2-\Phi)\\
\Phi_t=x-(N-2)\Phi+bme^{(2-q)t}|\Phi|^q.
\EA\right.$$
Setting $(2-q)a=\ln b$ , $\gt=t+a$, $x^{(a)}(\gt)=x(t-a)$ and $\Phi^{(a)}(\gt)=\Phi(t-a)$, we obtain that 
$$\left\{\BA{lll}
x^{(a)}_\gt(\gt)=x^{(a)}(\gt)(2-\Phi^{(a)}(\gt))\\
\Phi^{(a)}_\gt(\gt)=x^{(a)}(\gt)-(N-2)\Phi^{(a)}(\gt)+me^{(2-q)\gt}|\Phi^{(a)}(\gt)|^q.
\EA\right.$$
Equivalently, the function $\gr\mapsto U_1^{(a)}(\gr):=\ln\left(\gr^{-2}x^{(a)}(\ln\gr)\right)$ is a radial solution of $(\ref{Na-1})$ satisfying $(\ref{Na-23})$. We  show also that the uniqueness of the trajectory implies the uniqueness of $u$. Indeed, if $\hat u=U_2$ is another solution of $\ref{Na-1}$ with corresponding solution $(\hat x,\hat\Phi,\hat\Gth)$ of $(\ref{S20})$ where we have $\hat\Gth(t)=\frac1be^{(2-q)t}=e^{(2-q)t(t-a)}=\Gth(t-a)$, then $(\hat x,\hat\Phi,\hat\Gth)$ converges to $P_0$ when $t\to-\infty$ and admits the trajectory $\CT_1$, therefore there exists $h\in\BBR$ such that 
$(\hat x(t),\hat\Phi(t),\hat\Gth(t))=(x(t+h),\Phi(t+h),\Gth(t+h))$. Thus $h=-a$ and $(\hat x(t),\hat\Phi(t)=(x(t-a),\Phi(t-a))=(x^{(a)}(t),\Phi^{(a)}(t))$, and by construction $\hat u=u^{(a)}$.\smallskip

\subsubsection { Proof of Assertion 2} The proof is similar to the one of \cite[Theorem 17]{BV-V25} where the sign of $m$ does not matter. We recall the main streams of this proof. Setting $U=-u$, we look for a radial positive function $U$ satisfying for some $\gg<0$ and $\gr>0$,
\bel{U2}\left\{\BA{lll}\dsps
-U_{rr}-\frac{N-1}rU_r+m|U_r|^q+e^{-U}=0\quad\text{in }(0,\gr)\\[2mm]
\phantom{-}\dsps \lim_{r\to 0}r^{N-2}U(r)=-\gg\\[2mm]
\phantom{-----}U(\gr)=0,
\EA\right.\ee
and satisfying furthermore $\dsps \lim_{r\to 0}r^{N-1}U(r)=(N-1)\gg$. Then $V(r)=U_r(r)$ is expressed by
\bel{U3}
V(r)=(N-1)\gg r^{1-N}+ r^{1-N}\int_0^r\left(m|V|^q+e^{-U}\right)s^{N-1}ds,
\ee
where
\bel{U4}
U(r)=-\int_r^\gr V(s)ds.
\ee
We define the operator $(U,V)\mapsto \CK(U,V):=(\CK_1(U,V),\CK_2(U,V))$ with 
\bel{U5}\left\{\BA{lll}\dsps
\CK_1(U,V)(r)=-\int_r^\gr V(s)ds\\[2mm]\dsps
\CK_2(U,V)(r)=(N-1)\gg r^{1-N}+\int_0^r\left(m|V|^q+e^{-|U|}\right)s^{N-1}ds,
\EA\right.\ee
on the subspace $\CH$ of $C((0,\gr])\ti C((0,\gr])$ endowed with the norm
$$\norm{(U,V)}_\CH=\max\left\{\gs\sup_{0<r\leq \gr}r^{N-2}|U(r)|,\sup_{0<r\leq \gr}r^{N-1}|V(r)|\right\}:=\max\{\gs N_1(U),N_2(V)\}
$$
for some $0<\gs<1$ to be specified.\\
 The first estimate is a {\it Lipschitz estimate} essentially proved in \cite[formula (59)]{BV-V25}:
\bel{U6}\BA{lll}\dsps
\norm{\CK(U_1,V_1)-\CK(U_2,V_2)}_\CH\leq \max\left\{\frac{\gs}{N-2}N_2(V_1-V_2),\right.\\[2mm]\dsps
\phantom{------}
\left.\frac{\gr^2}{2}N_1(U_1-U_2)+\frac{mq\max\{N_2^{q-1}(V_1),N_2^{q-1}(V_2)\}}{N-q(N-1)}N_2(V_1-V_2)\right\}.
\EA\ee
 The second estimate is a {\it bound estimate}, see \cite[formula (63)]{BV-V25}. Observing that if $\norm{(U,V)}_\CH\leq R$ we have
 $|U(r)|\leq  R\gs^{-1}r^{2-N}$ and $|V(r)|\leq Rr^{1-N}$ which proves  that
\bel{U7}\norm{K(U,V)}_\CH\leq\max\left\{\frac{\gs R}{N-2},\frac{\gr^2R}{2\gs}+(N-2)|\gg|+\frac{\gr^N}{N}+\frac{mR^q\gr^{N-q(N-1)}}{N-q(N-1)}
\right\}.\ee
 Therefore, taking $\gs=\frac 34$, for any $R>0$ there exists $0<\gr_0<1$ and $k_0>0$ such that for any $|\gg|\leq k_0$ and $0<\gr\leq\gr_0$, there holds
 \bel{U8}
 \norm{(U,V)}_\CH\leq R\Longrightarrow \norm{K(U,V)}_\CH\leq R \ee

\nind With this Lipschitz estimate we can conclude the proof as follows:\smallskip

 \nind If $(\ref{U8})$ holds and $\norm{(U_1,V_1)}_\CH\leq R$ and $\norm{(U_2,V_2)}_\CH\leq R$, we have
 \bel{U9}\BA{lll}\dsps
 \norm{\CK((U_1,V_1)-\CK((U_2,V_2)}_\CH\leq \max\left\{\frac 34N_2(V_1-V_2),\right.\\[4mm]
 \phantom{---------------}\dsps\left.\frac{\gr^2}{2}N_1(U_1-U_2)+\frac{mqR^{q-1}}{N-q(N-1)}N_2(V_1-V_2)\right\}.
 \EA\ee
 Up to reducing $R$ up to $R\leq 1$ we can assume that $\frac{mqR^{q-1}}{N-q(N-1)}\leq \frac 14$ and obtain
 $$\max\left\{\frac 34N_2(V_1-V_2),\frac{\gr^2}{2}N_1(U_1-U_2)+\frac 14N_2(V_1-V_2),\right\}\leq \max\left\{\frac 34,\gr^2\right\}\norm{U-V}_\CH.
 $$
 Since $\gr<1$ the mapping $\CK$ is a contraction and it admits a fixed point $(U,V)$ such that $\norm{U,V}_\CH\leq R$. Then $U$ solves
 $$-U_{rr}-\frac{N-1}{r}U_r+m|U_r|^q+e^{-|U|}=0\quad\text{in }(0,\gr),
 $$
 and satisfies $\lim_{r\to 0}r^{N-2}U(r)=-\gg$ and $U(\gr)=0$. Clearly $U$ is monotone, therefore it is positive on $(0,\gr)$ and thus it satisfies $(\ref{U2})$.
 
 \subsubsection{Proof of Assertion 3}
We use again $U=-u$. For $\tilde m>m$ and $A\in\BBR$ we consider the function
 \bel{X1}
V_{A,\tilde m}(r)=\frac{q-1}{q-2}\left(-\frac{\gk}{\tilde m}\right)^{\frac{1}{q-1}}\left(r^{-\gb}-r_0^{-\gb}\right)+a
\ee
It satisfies
 \bel{X2}\BA{lll}
-\Gd V_{A,\tilde m}+m|\nabla V_{A,\tilde m}|^q+e^{-V_{A,\tilde m}}=(m-\tilde m)|\nabla V_{A,\tilde m}|^q+e^{-V_{A,\tilde m}}\quad\text{in }B_{r_0}\setminus\{0\}\\
\phantom{-\Gd |\nabla V_{A,\tilde m}|^qa+e^{-V_{A,\tilde m}}}
V_{A,\tilde m}(r_0)=A.
\EA\ee
We have that
$$(m-\tilde m)|\nabla V_{A,\tilde m}(r)|^q=(m-\tilde m)\left(-\frac{\gk}{\tilde m}\right)^{\frac{q}{q-1}}r^{-\frac q{q-1}}.
$$
Since
$$e^{-V_{A,\tilde m}(r)}=e^{-A}e^{\frac{q-1}{q-2}\left(-\frac{\gk}{\tilde m}\right)^{\frac{1}{q-1}}\left(r_0^{-\gb}-r^{-\gb}\right)}
$$
if we take
  \bel{X3}\BA{lll}\dsps
(\tilde m-m)\left(-\frac{\gk}{\tilde m}\right)^{\frac{1}{q-1}}r_0^{-\frac{q}{q-1}}\geq e^{-A}
\EA\ee
we have that 
 \bel{X4}\BA{lll}
-\Gd V_{A,\tilde m}+m|\nabla V_{A,\tilde m}|^q+e^{-V_{A,\tilde m}}\leq 0\quad\text{in }B_{r_0}\setminus\{0\}.
\EA\ee
On the other hand
 \bel{X5}\BA{lll}
-\Gd V_{A,m}+m|\nabla V_{A, m}|^q+e^{-V_{A, m}}=e^{-V_{A, m}}\geq 0\quad\text{in }B_{r_0}\setminus\{0\}.
\EA\ee
Since $V_{A,m}\geq V_{A,\tilde m}$ it follows by standard truncation of the domain techniques \cite[Theorem 1.4.6]{} that there exists a function $U_a$ (radial as $V_{A, m}$ and $V_{A, \tilde m}$) satisfying $(\ref{S28*})$ in $B_{r_0}\setminus\{0\}$ and 
  \bel{X6}\BA{lll}
V_{A, \tilde m}\leq U\leq V_{A, m}\quad\text{in }B_{r_0}\setminus\{0\}.
\EA\ee
and clearly $U(x)=r_0$ if $|x|=a$. Therefore $u_A=-U_A$ satisfies $(\ref{Na-1})$ in $B_{r_0}\setminus\{0\}$ and
  \bel{X6}\BA{lll}
-V_{A,m}\leq u_A\leq -V_{A, \tilde m}.
\EA\ee
By \rth{supth6-ter}-(3-1) we conclude that $u_A(r_0)=-A:=B$ satisfies $(\ref{Na-17-4})$. From this fact, the estimate of the gradient is standard.
 \qeda
  
\subsection{Proof of \rth{supth10}}
The existence of a radial solution satisfying $(\ref{Na-25})$ is an adaptation of  \cite[Theorem 28]{BV-V25}$
$ in which only the sign of $m$ is different ($m<0$ there), which has no impact on the construction. We define the functions
\bel{W24}
Z=-\frac{re^{u}}{u_r},V=r\left\vert u_r\right\vert ^{q-1}%
,\Phi=-ru_r,%
\ee
and look for a solution verifying $Z>0$ and $\Phi>0.$. Then $(Z,V,\Phi)$ verifies a
{\it quadratic system of order 3}
\bel{W25}\left\{
\begin{array}
{lll}
Z_{t}=Z(N-\Phi-mV-Z)\\
V_{t}=V(N-(N-1)q+(q-1)(mV+Z))\\
\Phi_{t}=\Phi(2-N+mV+Z)
\end{array}\right.
\ee
It is exactly the same system as in \cite[Theorem 28]{BV-V25}, the difference is that in that case  we were interested in solutions $u$
increasing, such that $\Phi<0,$ and $Z<0$. The system $(\ref{W25})$ admits a stationary point
$P_{0}=(0,V_{0},0)$ where $V_{0}=\frac{(N-1)q-N}{m(q-1)}=\frac{\kappa}{m}.$ 
Setting $V=V_{0}+\overline{V},$ the  linearised system at $P_{0}$ takes the
form%
\bel{W25bis}\left\{
\begin{array}
{lll}
Z_{t}=\frac{q}{q-1}Z\\
\overline{V}_{t}=(q-1)\kappa(\overline{V}+\frac{Z}{m})\\
\Phi_{t}=\frac{q-2}{q-1}\Phi
\end{array}\right.
\ee
The eigenvalues of the linear system $(\ref{W25bis})$ are 
$$
\lambda_{1}=\frac{q}{q-1},\quad\lambda_{2}=(N-1)q-N,\quad\lambda_{3}%
=\frac{q-2}{q-1}.%
$$
Since  $q>2,$ there holds $~0<\lambda_{3}<1<\lambda_{1}<2,$ and $\lambda
_{2}>N-2.$, and when $N\geq4$ they are all distincts. Hence  there exists an
infinity of trajectories such that $\Phi>0,$ and $Z>0.$ The eigenvalues
$\lambda_{1},\lambda_{2}$ are distinct, except in the case $N=3,$
$q=\frac{3+\sqrt{3}}{2},$ or $N=2,q=2+\sqrt{2}$. Since all the eigenvalues are
positive, there exists a neighborhood $\mathcal{V}$ of $P_{0}$ such that all
the trajectories starting from $\mathcal{V}$ converge to $P_{0}$ as
$t\rightarrow\infty.$ Adapting carefully the proof of \cite[Theorem 28]{BV-V25}, since the system of order 3 is not equivalent to the
equation $(\ref{Na-1})$ we show that {\it for any}
$u_{0}$ there exist a solution $u$ satisfying $(\ref{Na-25})$.\qeda

\medskip

\nind \Remark By a technical adaptation of the above proof it is possible to extend the result to the case $N=1$. The only difference is that the constant 
$c_{N,q,m}=\frac{q-1}{q-2}\left(\frac\gk m\right)^{\frac 1{q-1}}$ in $(\ref{Na-25})$ has to be replaced by $c_{1,q,m}=\frac{1-q}{q-2}\left(-\frac\gk m\right)^{\frac1{q-1}}<0$.

\section{Appendix}
In this appendix we give an alternative proof of a variant of Assertion 2 in \rth{supth8}, which does not use the notion of super and sub viscosity solutions but is based upon the geometric theory of dynamical systems. This proof is longer but its field of applications is not restricted to equation $(\ref{Na-1})$  and is well adapted to dynamical systems in presence of a saddle point equilibrium. The proof is inspired by a 
construction due to Dieudonn\'e \cite{Dieu}.
 \bth{Appendix} Let $N\geq 3$ and $q>2$. We can find $r_0>0$ and $\ge_0>0$ such that for any $|u_0|<\ge_0$ there exists a unique radial solution of $(\ref{Na-1})$ in $B_{r_0}^c$ satisfying $u(r_0)=u_0$ and 
    \bel{App1} \lim_{r\to\infty}r^{\gb}u(r)=-\frac {q-1}{q-2}\left(\frac \gk m\right)^{\frac 1{q-1}}\quad\text{and }\,\lim_{r\to\infty}r^{\frac 1{q-1}}u_r(r)=-\left(\frac \gk m\right)^{\frac 1{q-1}}.\ee 
\es
\Proof 1- {\it The dynamical system formulation.} We use system $(\ref{T3})$ with variable $U(r)=-u(r)$, $\xi(t)=r^{-\gb}U(r)$, $\eta(t)=r^{\frac{1}{q-1}}U_r(r)$ and $t=\ln r$. We define the vector field
$$\CF(\xi,\eta)=(\gb\xi+\eta,-\gk\eta+m|\eta|^q).
$$
with equilibria $(0,0)$ and  $(\xi_m,\eta_m):=(\xi_m,\frac{q-2}{q-1}\xi_m)$ where $\xi_m=\frac{q-1}{q-2}\left(\frac\gk m\right)^{\frac1{q-1}}$. Then
\bel{W1}\BA{lll}\dsps 
D\CF(\xi_m,\eta_m)=\begin{pmatrix}\gb&1\\[2mm]
0&-\gk+qm\eta_m^{q-1}\end{pmatrix}=\begin{pmatrix}\gb&1\\[2mm]
0&(q-1)\gk
\end{pmatrix}:=Q,\EA\ee
and the eigenvalues of $Q$ are $\gb<0<(q-1)\gk$ and corresponding eigenvectors $(1,0)$ and $((q-1)\gk-\gb,1)$. We have seen that $(\xi_m,\eta_m)$ is a saddle point. We set $x(t)=\xi(t)-\xi_m$, $y(t)=\eta(t)-\eta_m$ and obtain the system
\bel{W2}\left\{\BA{lll}\dsps 
x_t=\gb x+y\\[1mm]
y_t=(q-1)\gk y+b_1(y)+b_2(t,x(t)),
\EA\right.\ee
with
\bel{W3}
b_1(y)=\frac{q(q-1)}{2}(\gk^{q-2}m)^{\frac{1}{q-1}}+o(y^2)\,\text{ as }\,y\to 0\,\text{ and }\,b_2(t,x)=e^{\frac{qt}{q-1}}e^{(-\xi_m+x)e^{-\gb t}},
\ee
The matrix $Q$ of the linear part of the system $(\ref{W2})$ can be put into a diagonal form. If we set $h=\frac{1}{(q-1)\gk-\gb}>0$ and 
\bel{W3-0}
\left\{\BA{lll} &x=X+Y\\ &y=\frac{Y}{h}\EA\right.\;\text{ which is equivalent to }\; \left\{\BA{lll} &X=x-hy\\ &Y=hy.\EA\right.\ee
then
\bel{W3-1}
\left\{\BA{lll}
X_t=\gb X-h\left(b_1(h^{-1}Y)+ b_2(t,X+Y)\right)\\
Y_t=(q-1)\gk Y+h\left(b_1(h^{-1}Y)+ b_2(t,X+Y)\right).
\EA\right.
\ee
We set
$$B(t,X,Y)=b_1(h^{-1}Y)+ b_2(t,X+Y).
$$
Since by $(\ref{W3})$ we have
$$|b_1(h^{-1}Y)|=o(Y^2)\quad\text{as }Y\to 0
$$
we deduce that for any $\ge>0$ there exists $\gd>0$ such that for any $(Y,Y')\in\BBR^2$ satisfying $\max\{|Y|, |Y'|\}\leq \gd$, we have 
\bel{W6}
|b_1(h^{-1}Y)-b_1(h^{-1}Y')|\leq \ge(|Y-Y'|).
\ee
Next we estimate the second term in $B(t,X,Y)$.
\blemma{Lem2} For any $M>0$ there exists $a>0$ such that for any $(X,X',Y,Y')\in\BBR^4$ satisfying $\max\{|X|,|X'|, |Y|, |Y'|\}\leq M$, we have 
\bel{W7}\BA{lll}
|b_2(X+Y,t)-b_2(X'+Y',t)|\leq ae^{\frac {qt}{q-1}-\xi_me^{-\gb t}}\left(|X-X'|+|Y-Y'|\right).
\EA\ee
\es

\nind We set $\CU(t)=(X(t),Y(t))$ and
\bel{W8}
F(t,X(t),Y(t))=F(t,\CU(t))=(-hB(t,\CU(t),hB(t,\CU(t))).
\ee
The main difference (and difficulty) with Dieudonn\'e's approach \cite[Chapter 13]{Dieu} is that $F(t,0)$ is not zero. We look for solutions satisfying the following integral equation on $[t_0,\infty)$ for some $t_0$ to be fixed later on,
\bel{W9}\left\{\BA{lll}\dsps 
X(t)=e^{\gb (t-t_0)}x_0-h\int_{t_0}^te^{\gb (t-s)}B(s,X(s),Y(s))ds\\[3mm]
\dsps
Y(t)=-h\int_t^\infty e^{(q-1)\gk(t-s)}B(s,X(s),Y(s))ds,
\EA\right.\ee
for some $x_0$ and provided $X^2(t)+Y^2(t)$ is small enough, and more precisely
\bel{W10}\left\{\BA{lll}\dsps 
|x_0|\leq \frac{\xi_m}{3}\\[2mm]
\dsps
|X(t)|+|Y(t)|\leq 2\sqrt{X^2(t)+Y^2(t)}\leq \frac{\xi_m}{3}.
\EA\right.\ee

\nind 2- {\it  The iterative scheme .} We set $\CU_n(t)=(X_n(t),Y_n(t))$ and consider the scheme associated to the integral equations defined in $(\ref{W9})$   where $\CU_0=(0,0)$ and 
  \bel{W11}\left\{\BA{lll}\dsps 
X_{n+1}(t)=e^{\gb (t-t_0)}x_0-h\int_{t_0}^te^{\gb (t-s)}B(s,X_n(s),Y_n(s))ds\\[3mm]\dsps
Y_{n+1}(t)=-h\int_{t_0}^te^{\gb (t-s)}B(s,X_n(s),Y_n(s))ds.
\EA\right.\ee
When $n=1$ we have for $t\geq t_0$
\bel{W12}
|X_1(t)|\leq e^{\gb (t-t_0)}\left(|x_0|+e^{\gb t_0}h\int_{t_0}^te^{2s-\xi_me^{-\gb s}}ds\right)\leq e^{\gb (t-t_0)}\left(|x_0|+c_2(t_0)\right),
\ee
where 
$$c_2(t_0)=e^{\gb t_0}\int_{t_0}^\infty e^{2s-\xi_me^{-\gb s}}ds=\int_{e^{t_0}}^\infty e^{-\xi_me^{\gb X}}XdX\to 0\quad\text{as }t_0\to\infty,$$
and 
\bel{W13}
|Y_1(t)|\leq e^{\gk(q-1)t}\left(c_3h\int_{t}^\infty e^{\frac{q-2-\gk(q-1)^2}{q-1}s-\xi_me^{-\gb s}}ds\right)\leq c_4(t_0)e^{-\frac{\xi_m}{2}e^{-\gb t}},
\ee
where also $c_4(t_0)\to 0$ if $t_0\to\infty$.
\blemma{Lem3} There exists $K>0$ such that for any $n\geq 1$ there holds
\bel{W14}
|X_{n+1}(t)-X_{n}(t)|+|Y_{n+1}(t)-Y_{n}(t)|\leq \frac{K}{2^n}e^{\frac\gb2(t-t_0)}(|x_0|+c_2).
\ee
\es
\Proof For $n\geq 1$ we have
\bel{W15}\BA{lll}\dsps
|X_{n+1}(t)-X_{n}(t)|\leq h\int_{t_0}^te^{\gb(t-s)}\left(|b_1(h^{-1}Y_n(s))-b_1(h^{-1}Y_{n-1}(s))|^{\phantom{P^L}}\right.\\[4mm]
\phantom{------X_{n+1}(t)-X_{n}(t)} \dsps\left.+|b_2(s,X_n(s),Y_n(s))-b_2(s,X_{n-1}(s),Y_{n-1}(s))|^{\phantom{P^L}}\!\!\!\!\!\!\!\right)ds\\[4mm]
\phantom{|X_{n+1}(t)-X_{n}(t)| }\dsps
\leq \ge h\int_{t_0}^te^{\gb(t-s)}|Y_{n}(s)-Y_{n-1}(s)|ds\\[4mm]
\phantom{|X_{n+1}(t)-X_{n}(t)|+} \dsps
+ah\int_{t_0}^te^{\gb(t-s)}e^{\frac {qs}{q-1}-\xi_me^{-\gb s}}\left(|X_{n}(s)-X_{n-1}(s)|+|Y_n(s)-Y_{n-1}(s)|\right)ds
\\[4mm]
\phantom{|X_{n+1}(t)-X_{n}(t)|} \dsps
\leq A_1(t)+A_2(t).
\EA\ee 
By induction we assume that 
\bel{W16}
|X_{n}(t)-X_{n-1}(t)|+|Y_{n}(t)-Y_{n-1}(t)|\leq \frac{K}{2^{n-1}}e^{\frac\gb2(t-t_0)}(|x_0|+c_2).
\ee
Then
\bel{W17}\BA{lll}\dsps
A_1(t)\leq \ge h\int_{t_0}^te^{\gb (t-s)}\left(|X_{n}(s)-X_{n-1}(s)|+|Y_{n}(s)-Y_{n-1}(s)|\right)ds\\[4mm]
\phantom{I(t)}\dsps
\leq \frac{(q-1)K \ge h}{2^{n-2}(q-2)}e^{\gb(t-{t_0})}\left(|x_0|+c_2\right).
\EA\ee
We first choose $\ge h\leq \frac{q-2}{16(q-1)}$ and we obtain
\bel{W18}
A_1(t)\leq \frac K{2^{n+2}}e^{\gb(t-{t_0})}\left(|x_0|+c_2\right).
\ee
Concerning the term $A_2(t)$ we have
$$\BA{lll}A_2(t)\dsps \leq ahe^{\gb t_0}\int_{t_0}^te^{2s-\xi_me^{-\gb s}}\left(|X_{n}(s)-X_{n-1}(s)|+|Y_n(s)-Y_{n-1}(s)|\right)ds\\[4mm]
\phantom{B(t)}\dsps \leq \frac{Kah(|x_0|+c_2)}{2^{n-1}}e^{\frac\gb2(t-t_0)}\int_{t_0}^te^{\frac\gb2(t-s)+\frac{qs}{q-1}-\xi_me^{-\gb s}}ds.
\EA$$
Since for all $t_0\leq s\leq t$ there holds
$$e^{\frac\gb 2(t-s)+\frac{qs}{q-1}-\xi_me^{\frac{q-2}{q-1}s}}\leq e^{\frac\gb 2(t-s)+\frac{qs}{q-1}-\xi_me^{\frac{q-2}{q-1}t_0}},
$$
we can take $t_0$ large enough such that 
\bel{W19}\BA{lll}\dsps
A_2s(t)\leq \frac{K}{2^{n+2}}e^{\frac\gb 2(t-{t_0})}\left(|x_0|+c_2\right).
\EA\ee
Next we estimate $|Y_{n+1}(t)-Y_{n}(t)|$. We have

\bel{W20}\BA{lll}\dsps
|Y_{n+1}(t)-Y_n(t)|\leq h\int_t^\infty e^{(q-1)\gk(t-s)}\left(|b_1(h^{-1}Y_n(s))-b_1(h^{-1}Y_{n-1}(s))|^{\phantom{P^L}}\right.\\[4mm]
\phantom{-------Y_{n+1}(t)-Y_{n}(t)} \dsps\left.+|b_2(s,X_n(s),Y_n(s))-b_2(s,X_{n-1}(s),Y_{n-1}(s))\right)ds\\[4mm]
\phantom{|Y_{n+1}(t)-Y_n(t)|}\dsps
\leq \ge h\int_t^\infty e^{(q-1)\gk(t-s)}|Y_{n}(s)-Y_{n-1}(s)|ds\\[4mm]
\phantom{|_{n}(t)|} \dsps
+ah\int_t^\infty e^{(q-1)\gk(t-s)}e^{\frac {qs}{q-1}-\xi_me^{-\gb s}}\left(|X_{n}(s)-X_{n-1}(s)|+|Y_n(s)-Y_{n-1}(s)|\right)ds\\[4mm]
\phantom{|Y_{n+1}(t)-Y_n(t)|} \dsps
\leq A_3(t)+A_4(t).
\EA\ee
By the induction assumption $(\ref{W16})$,
$$\BA{lll}\dsps A_3(t)\leq \frac{K \ge h }{2^{n-1}}(|x_0|+c_2)\int_t^\infty e^{(q-1)\gk(t-s)}e^{\frac\gb 2(s-{t_0})}ds\\[4mm]
\phantom{C(t)}\dsps 
\leq \frac{(q-1)K \ge h}{2^{n-2}(2\gk(q-1)^2+q-2)}(|x_0|+c_2)e^{\frac\gb 2(t-t_0)}.
\EA$$
We next choose $\ge h\leq \min\left\{\frac{q-2}{16(q-1)},\frac{2\gk(q-1)^2+q-2}{16(q-1)}\right\}$ and we obtain
\bel{W21}\BA{lll}\dsps
A_3(t)\leq \frac{K}{2^{n+2}}e^{\frac\gb 2(t-{t_0})}\left(|x_0|+c_2\right).
\EA\ee
For the last term, we write
\bel{W22}\BA{lll}\dsps
A_4(t)\leq \frac{Ka h}{2^{n-1}}\left(|x_0|+c_2\right)\int_t^\infty e^{(q-1)\gk(t-s)+\frac {qs}{q-1}-\xi_me^{-\gb s}+\frac\gb 2(s-t_0)}ds\\[4mm]
\phantom{D(t)}\dsps
\leq \frac{K h}{2^{n+2}}e^{\frac\gb 2(t-{t_0})}\left(|x_0|+c_2\right)8a
\int_t^\infty e^{((q-1)\gk-\gb)(t-s)+\frac {qs}{q-1}-\xi_me^{-\gb s}}ds
\EA\ee
Since
$$\lim_{s\to\infty}\left(\frac {qs}{q-1}-\xi_me^{-\gb s}\right)=-\infty
$$
it follows that 
$$\lim_{t\to\infty}\int_t^\infty e^{((q-1)\gk-\gb)(t-s)+\frac {qs}{q-1}-\xi_me^{-\gb s}}ds=0
$$
Hence there exists $t_0$ such that for all $t\geq t_0$,
$$\int_t^\infty e^{((q-1)\gk-\gb )(t-s)+\frac {qs}{q-1}-\xi_me^{-\gb s}}ds\leq \frac{1}{8a}.
$$
This implies
\bel{W23}\BA{lll}\dsps
D(t)\leq \frac{K}{2^{n+2}}e^{\frac\gb 2(t-{t_0})}\left(|x_0|+c_2\right).
\EA\ee
Combining $(\ref{W18})$, $(\ref{W19})$, $(\ref{W21})$ and $(\ref{W23})$ we obtain
$(\ref{W14})$.\qeda\smallskip

Next we define the weighted space $C_w([t_0,\infty))$ with $w(t)=e^{-\frac\gb 2(t-t_0)}$ by 
$$\dsps C_w(t_0,\infty)=\left\{(\phi,\psi)\in C([t_0,\infty))\ti C([t_0,\infty)) :\sup_{t\geq t_0}w(t)(|\phi(t)|+|\psi(t)|)<\infty\right\},
$$
and we denote by $\norm {(\phi,\psi)}_{C_w(t_0,\infty)}$ this last quantity.
\bprop{conv} The sequence $\{(X_n,Y_n)\}$ is a Cauchy sequence in $C_w([t_0,\infty))$. Its limit $\CU(t)=(X(t),Y(t))$ satisfies $(\ref{W9})$. Hence 
$x=X+Y$ and $y=\frac Yh$
satisfy $(\ref{W2})$ and $\dsps\lim_{t\to\infty}(x(t),y(t))=(0,0)$. Thus $(\xi(t),\eta(t))$ is a solution of system $(\ref{T3})$ on $(t_0,\infty)$ such that $\dsps\lim_{t\to\infty}(\xi(t),\eta(t))=(\xi_m,\eta_m)$.
\es
\Proof The space $C_w([t_0,\infty))$ is complete and we have
$$\norm {(X_{n+p}-X_n,Y_{n+p}-Y_n)}_{C_w(t_0,\infty)}\leq \frac 1{2^{n-1}}\left(|x_0|+c_2\right).
$$
Hence there exists $(X(t),Y(t))=\lim_{t\to\infty}(X_n(t),Y_n(t))$ in this space $C_w([t_0,\infty))$. By $(\ref{W6})$ and \rlemma {Lem2}, we have 
$\lim_{n\to\infty} F(t,\CU_n(t))=F(t,\CU(t))$ in $C_w([t_0,\infty))$. Therefore $\CU(t)=(X(t),Y(t))$ satisfies the integral equation $(\ref{W9})$ which is equivalent to the system $(\ref{W3-1})$. By the transformation  $(\ref{W3-0})$, $(x,y)$ satisfies 
$(\ref{W2})$ and $(x(t),y(t))\to (0,0)$ when $t\to\infty$.\qeda\\

Returning to the variable $u$ by the transformation defining $\xi$ and $\eta$, the proof of Assertion 2 in \rth{supth8} follows with $r_0=e^{t_0}$.\qeda\\

\nind\Remark The above method can be adapted to prove a slightly different statement than Assertion 3 in \rth {supth9}.\medskip

\begin{center}
{\bf Acknowledgement}
\end{center}
The second author is supported by Fondecyt grant 1250522 from ANID, Chile. \medskip

\begin{center}

 {\bf STATEMENTS AND DECLARATION}\smallskip
 
 \end{center}
 \nind {\bf Conflict of interest}: The authors declare that they have no conflict of interest and that this article follows all the ethical rules.\smallskip

\nind {\bf Data availability}: Data sharing is not applicable to this article as no datasets were generated or analysed during the study.

\nind Marie-Fran\c{c}oise Bidaut-V\'eron \\
Institut Denis Poisson CNRS UMR 7013. \\
Universit\'e de Tours\\
Tours, France\\
{\it veronmf@univ-tours.fr}\\

\nind Marta Garcia-Huidobro \\
Departamento de Matematicas\\ 
Pontifica Universidad Catolica de Chile\\
Santiago de Chile\\
{\it mgarcia@mat.puc.cl}\\

\nind Laurent V\'eron \\
Institut Denis Poisson CNRS UMR 7013\\
Universit\'e de Tours\\
Tours, France\\
{\it veronl@univ-tours.fr}

 \end{document}